\documentclass[11pt,a4paper,leqno,twoside]{amsart}
\usepackage{latexsym,amssymb,amsmath}
\input xy
\xyoption{all}

\def\CC{\mathbb C}

\def\RR{\mathbb R}
\def\HH{\mathbb H}
\def\AA{{\mathbb A}}

\def\OO{\mathbb O}

\def\ZZ{\mathbb Z}

\def\11{\mathbf 1}
\def\PP{\mathbb P}

\def\e1{\varepsilon_1}
\def\e2{\varepsilon_2}
\def\e3{\varepsilon_3}

\def\P2{{\PP}^2}

\def\00{\underline{0}}
\def\J0{{\cal J}_3(\underline{0})}

\def\PJ0{\PP({\cal J}_3(\underline{0}))}

\def\a{\alpha}

\def\b{\beta}
\def\g{\gamma}
\def\s{\sigma}

\def\d{\delta}

\def\e{\varepsilon}

\def\AP2{{\AA\PP}^2}
\def\RP2{{\RR\PP}^2}
\def\CP2{{\CC\PP}^2}
\def\HP2{{\HH\PP}^2}
\def\OP2{{\OO\PP}^2}

\newtheorem{theo}{Theorem}[section]
\newtheorem{coro}[theo]{Corollary}
\newtheorem{lemm}[theo]{Lemma}
\newtheorem{prop}[theo]{Proposition}
\theoremstyle{remark}
\newtheorem{rema}[theo]{Remark}
\theoremstyle{definition}
\newtheorem{defn}{Definition}

\begin{document}

\title{Primarily quasilocal fields and 1-dimensional abstract local
class field theory}

\author
{I.D. Chipchakov}
\address{Institute of Mathematics and Informatics\\Bulgarian Academy
of Sciences\\Acad. G. Bonchev Str., bl. 8\\1113, Sofia, Bulgaria\\
e-mail chipchak@math.bas.bg} 
\thanks{The author was partially supported by Grant MI-1503/2005 of
the Bulgarian NSF}
\begin{abstract}
Let $E$ be a field satisfying the following conditions: (i) the
$p$-component of the Brauer group Br$(E)$ is nontrivial whenever $p$
is a prime number for which $E$ is properly included in its maximal
$p$-extension; (ii) the relative Brauer group Br$(L/E)$ equals the maximal
subgroup of Br$(E)$ of exponent $p$, for every cyclic extension
$L/E$ of degree $p$. The paper proves that finite abelian
extensions of $E$ are uniquely determined by their norm groups and
related essentially as in the classical local class field theory.
This includes analogues to the fundamental correspondence, the local
reciprocity law and the local Hasse symbol.
\end{abstract}

\maketitle

\section{Introduction and statements of the main results}

\medskip
This paper is concerned with finite abelian extensions of primarily
quasilocal (abbr, PQL) fields, and can be viewed as a continuation
of \cite{Ch6}, I. When $E$ is a strictly PQL-field, it shows that
these extensions and their norm groups are related as in the
fundamental correspondence of the classical local class field theory
(see \cite{I}, page 101). The paper proves that they are subject to
an exact analogue to the local reciprocity law (formulated, e.g., in
\cite{Sch}, Ch. 6, Theorem~8, and \cite{I}, Theorem~7.1), and to a
partial analogue to the local Hasse symbol (as characterized in
\cite{Ko2}, Ch. 2, see also \cite{Wh1}, and \cite{I}, Theorems~6.9
and 6.10) of form determined by invariants of the Brauer group
Br$(E)$. It takes a step towards characterizing the fields whose
finite abelian extensions have the above properties. When $E$
belongs to some special classes of traditional interest, the present
research enables one to achieve this aim and to find a fairly
complete description of the norm groups of arbitrary finite
separable extensions of $E$ (see Section 3, \cite{Ch5,Ch9} and the
references in \cite{Ch10}, Remark~3.9).

\medskip
The basic notation, terminology and conventions kept in this paper
are standard and virtually the same as in \cite{Ch6}, I, and
\cite{Ch10}. Throughout, $\mathbb P$ denotes the set of prime
numbers and every algebra $A$ is understood to be associative with a
unit lying in the considered subalgebras of $A$. Simple algebras are
supposed to be finite-dimensional over their centres, Brauer and
value groups are written additively, Galois groups are viewed as
profinite with respect to the Krull topology, and the considered
profinite group homomorphisms are continuous. For each field $E$, $E
^{\ast }$ denotes its multiplicative group, $E _{\rm sep}$ a
separable closure of $E$, $\mathcal{G} _{E} = \mathcal{G}(E _{\rm
sep}/E)$ is the absolute Galois group of $E$, $s(E)$ stands for the
class of central simple $E$-algebras, $[B]$ denotes the similarity
class of any $B \in s(E)$, $d(E)$ is the subclass of division
algebras $D \in s(E)$, and $P(E) = \{p \in \mathbb P\colon \ E(p)
\neq E\}$, where $E(p)$ is the maximal $p$-extension of $E$ in $E
_{\rm sep}$. As usual, $E$ is said to be formally real, if $-1$ is
not presentable over $E$ as a finite sum of squares; $E$ is called
nonreal, otherwise. We say that $E$ is Pythagorean, if it is
formally real and the set $E ^{\ast 2} = \{\lambda ^{2}\colon \ 
\lambda \in E ^{\ast }\}$ is additively closed. For any field 
extension $F/E$, $\rho _{E/F}$ denotes the scalar extension map of 
Br$(E)$ into Br$(F)$, Br$(F/E)$ the relative Brauer group of 
$F/E$, and $I(F/E)$ is the set of intermediate fields of $F/E$. When 
$F/E$ is finite and separable, Cor$_{F/E}$ stands for 
the corestriction map of Br$(F)$ into Br$(E)$ (see \cite{T}). We say
that $E$ is $p$-quasilocal, for some $p \in \mathbb P$, if it
satisfies one of the following conditions: (i) every cyclic degree
$p$ extension of $E$ is embeddable as a subalgebra in each $\Delta
\in d(E)$ of (Schur) index $p$; (ii) the $p$-component Br$(E) _{p}$
of Br$(E)$ is trivial or $p \notin P(E)$. The field $E$ is called
PQL, if it is $p$-quasilocal for every $p \in \mathbb P$; when this
holds, $E$ is said to be strictly PQL in case Br$(E) _{p} \neq
\{0\}$, $p \in P(E)$. We say that $E$ is quasilocal, if its finite
extensions are PQL.

\medskip
Local fields and $p$-adically closed fields are strictly quasilocal
(abbr, SQL), i.e. their finite extensions are strictly PQL (cf.
\cite{S1}, Ch. XIII, Sect. 3, \cite{PR}, Theorem~3.1 and Lemma~2.9,
and \cite{Ch4}, Sect. 3). The strictly PQL-property has been fully
characterized in the following two classes: (i) algebraic extensions
of global fields \cite{Ch9}; (ii) Henselian discrete valued fields
\cite{Ch4}, Sect. 2. These facts extend the arithmetic basis of this
research which is further motivated by results of \cite{Ch6}, I, on
Brauer groups of PQL-fields and on absolute Galois groups of
quasilocal fields. For reasons clarified in the sequel, our approach
to the main topic of this paper is purely algebraic. Our starting
point is the fact that residue fields of Henselian valued stable
fields are PQL in the case of totally indivisible value groups, in
the sense of (2.1) (i). Specifically, a Henselian discrete valued
field $(K, v)$ with a perfect residue field $\widehat K$ is stable
if and only if $\widehat K$ is stable and PQL \cite{Ri}, 
Proposition~2 (cf. also \cite{Ch6}, I, and \cite{Ch11}, 
Proposition~2.3). These and other related results show, with their 
proofs, that PQL-fields partly resemble local fields in a number of 
respects (see Propositions 2.2-2.3, \cite{Ch6}, I, Lemma~4.3 and 
\cite{Ch6}, II, Lemma~2.3). For example, by \cite{Ch6}, I, 
Lemma~4.3, if $\widehat K$ is $p$-quasilocal and $\widetilde L _{1}$ 
and $\widetilde L _{2}$ are different extensions of $\widehat K$ in 
$\widehat K (p)$ of degree $p$, then the inner product $N(\widetilde 
L _{1}/\widehat K)N(\widetilde L _{2}/\widehat K)$ of the norm 
groups $N(\widetilde L _{i}/\widehat K)$, $i = 2$, is equal to 
$\widehat K ^{\ast }$. This result and its key role in the proof of 
\cite{Ch6}, I, Theorem~4.1, attract interest in the study of the 
PQL-property along the lines of the classical local class field 
theory, with the notion of a local field extended as follows:

\medskip
\begin{defn}
Let $E$ be a field, Nr$(E)$ the set of norm groups of $E$, and
$\Omega (E)$ the set of finite abelian extensions of $E$ in $E _{\rm
sep}$. We say that $E$ admits $1$-dimensional local class field
theory (abbr, LCFT), if the natural mapping of $\Omega (E)$ into
Nr$(E)$ (by the rule $M \to N(M/E)$, $M \in \Omega (E)$)  is
injective and the following condition holds, for each $M _{1}$, $M
_{2} \in \Omega (E)$:

\medskip
(1.1) The norm group (over $E$) of the compositum $M _{1}M _{2}$
equals the intersection $N(M _{1}/E) \cap N(M _{2}/E)$, and $N(M
_{1} \cap M _{2}/E) = N(M _{1}/E)N(M _{2}/E)$.
\end{defn}

We say that $E$ is a field with $1$-dimensional local $p$-class
field theory (abbr, local $p$-CFT), for a given $p \in \mathbb P$,
if the fields from the set $\Omega _{p} (E) = \{L \in \Omega
(E)\colon \ L \subseteq E(p)\}$ are uniquely determined by their
norm groups and satisfy condition (1.1). When this is the case and
$p \in P(E)$, we have Br$(E) _{p} \neq \{0\}$ (see Proposition \ref{prop3.3}).
Observe that $E$ admits LCFT if and only if it admits local $p$-CFT,
for every $p \in P(E)$. This follows from Lemma \ref{lemm2.1} and shows that
PQL-fields with LCFT are strictly PQL.

\medskip
The main purpose of this paper is to shed light on the place of
strictly PQL-fields in LCFT by proving the following:

\medskip
\begin{theo}
\label{theo1.1}
Strictly {\rm PQL}-fields admit {\rm LCFT}. Conversely, a field $E$
admitting such a theory is strictly {\rm PQL}, if each $D \in d(E)$ of
prime exponent $p$ is similar to a tensor product of cyclic division
$E$-algebras of index $p$.
\end{theo}

\medskip
\begin{theo}
\label{theo1.2}
Let $E$ and $M$ be fields, such that $E$ is strictly {\rm PQL},
$P(E) \neq \phi $  and $M \in \Omega (E)$. For each $p \in P(E)$,
let $_{p} {\rm Br}(E) = \{b _{p} \in {\rm Br}(E)\colon \ pb _{p} =
0\}$, $I _{p}$ be a basis and $d(p)$ the dimension of $_{p} {\rm
Br}(E)$ as a vector space over the field $\mathbb F _{p}$ with $p$
elements, $\mathcal{G}(M/E) _{p}$ the Sylow $p$-subgroup of the
Galois group $\mathcal{G}(M/E)$, and $\mathcal{G}(M/E) _{p} ^{d(p)}$
the direct product, indexed by $I _{p}$, of isomorphic copies of
$\mathcal{G}(M/E) _{p}$. Then the direct product $\mathcal{G}(M/E)
^{{\rm Br}(E)} = \prod _{p \in P(E)} \mathcal{G}(M/E) _{p} ^{d(p)}$
and the quotient group $E ^{\ast }/N(M/E)$ are isomorphic.
\end{theo}

\medskip
Before stating our third main result, recall that a field $F$ is
Pythagorean with $F (2) = F(\sqrt{-1})$ if and only if it is
formally real and $2$-quasilocal \cite{Ch6}, I, Lemma~3.5. Note also
that, by \cite{Wh2}, Theorem~2, if $p \in P(F)$, then $F (p)$
contains as a subfield a $\mathbb Z _{p}$-extension $U _{p}$ of $F$
(i.e. $U _{p}/F$ is Galois with $\mathcal{G}(U _{p}/F)$ isomorphic
to the additive group $\mathbb Z _{p}$ of $p$-adic integers) unless
$p = 2$ and $F$ is Pythagorean. We retain notation as in Theorem
\ref{theo1.2}.

\medskip
\begin{theo}
\label{theo1.3}
Let $E$ be a strictly {\rm PQL}-field, such that $P(E) \neq \phi $,
and let $E _{\infty } \subseteq E _{\rm sep}$ be the compositum of
fields $E _{p}$, $p \in P(E)$, where $E _{p}/E$ is a $\mathbb Z
_{p}$-extension, if $p > 2$ or $E$ is nonreal, and $E _{2} = E(2)$
when $E$ is formally real. Then there exists a set $H _{E} = \{( \ ,
M/E)\colon \ E ^{\ast } \to \mathcal{G}(M/E) ^{{\rm Br}(E)}$, $M \in
\Omega (E)\}$ of surjective group homomorphisms with the following
properties:
\par
{\rm (i)} The kernel of $( \ , M/E)$ coincides with $N(M/E)$, for
each $M \in \Omega (E)$;
\par
{\rm (ii)} If $M \in \Omega (E)$ and $K$ is an intermediate field of
$M/E$, then $( \ , K/E)$ equals the composition $\pi _{M/K} \circ (
\ , M/E)$, where $\pi _{M/K}\colon {\mathcal G}(M/E) ^{{\rm Br}(E)} \to
\mathcal{G}(K/E) ^{{\rm Br}(E)}$ is the homomorphism mapping the $i
_{p}$-th component of
\par\noindent
${\mathcal G}(M/E) _{p} ^{d(p)}$ on the $i
_{p}$-th component of $\mathcal{G}(K/E) _{p} ^{d(p)}$ as the natural projection
$\mathcal{G}(M/E) _{p} \to \mathcal{G}(K/E) _{p}$, for each pair $(p, i
_{p}) \in P(E) \times I _{p}$;
\par
{\rm (iii)} The set $H _{E}$ is uniquely determined by the mappings
$( \ , \Gamma /E)$, where $\Gamma $ runs through the set of finite
extensions of $E$ in $E _{\infty }$ of primary degrees.
\end{theo}

\medskip
Theorems \ref{theo1.1} and \ref{theo1.2} show the strong influence of Br$(E)$ on a
number of algebraic properties of a PQL-field $E$. They are
obtained from similar results on finite abelian $p$-extensions of
$p$-quasilocal fields, stated as Theorem \ref{theo3.1}. This approach uses
several properties of $p$-quasilocal fields without generally
valid analogues for PQL-fields (see Propositions \ref{prop2.2} and \ref{prop2.3},
\cite{Ch6}, I, Corollary~8.5, and \cite{Ch10}, Proposition~6.3). As
shown in \cite{Ch10,Ch11} and \cite{Ch6}, II, Sect. 3, it enables one
to describe the isomorphism classes of Brauer groups of major types
of PQL-fields, and of the reduced parts of Brauer groups of
equicharacteristic Henselian valued absolutely stable fields with
totally indivisible value groups. Thus it turns out that usually
powerful methods of valuation theory are virtually inapplicable to
many PQL and most presently known quasilocal fields (see (2.4) and
(2.5) (iii), and compare (2.1) with (2.2) and Remark \ref{rema2.4}). At the
same time, the proofs in \cite{Ch10} show at crucial points that the
study of the PQL-property can effectively rely on constructive
methods based on properties (established in the 1990's, see e.g.
\cite{FSS}, Proposition~2.6, and \cite{SVdB}) of relative Brauer
groups of extensions obtained as transfers of function fields of
Brauer-Severi and other varieties. This makes it possible to apply
Theorem \ref{theo3.1} and other results about $p$-quasilocal fields to Brauer
groups of arbitrary fields, and so leads to a better understanding
of the relations between Galois groups and norm groups of finite
Galois extensions of quasilocal fields (see Remark \ref{rema4.3} and
\cite{Ch10}).

\medskip
When $E$ is a local field, the former assertion of Theorem \ref{theo1.1}
yields the fundamental correspondence of the classical local class
field theory. If $E$ is merely strictly PQL with $d(p) =
1$, for all $p \in P(E)$ (i.e. with Br$(E)$ embeddable in
$\mathbb Q/\mathbb Z$, the quotient group of the additive
group of rational numbers by the subgroup of integers), then Theorem
\ref{theo1.2} states that $E ^{\ast }/N(M/E) \cong \mathcal{G}(M/E)$. This holds,
for instance, in the following cases:

\medskip
(1.2) (i) $E$ is an algebraic strictly PQL-extension of a global
field $E _{0}$; when Br$(E) _{p} \neq \{0\}$, by \cite{Ch9},
Theorem~2.1, $p \in P(E)$ and Br$(E) _{p}$ is isomorphic to the
quasicyclic $p$-group $\mathbb Z (p ^{\infty })$ unless $p = 2$ and
$E$ is formally real (big families of such $E$ can be constructed by
applying \cite{Ch9}, Theorem~2.2).
\par
(ii) The triple ${\rm AT}(E) = (\mathcal{G} _{E}, \{\mathcal{G}
_{F}\colon \ F \in \Sigma \}, E _{\rm sep} ^{\ast })$, $\Sigma $
being the set of finite extensions of $E$ in $E _{\rm sep}$, is an
Artin-Tate class formation (see \cite{AT}, Ch. 14).

\medskip
It is known that (1.2) (ii) holds, if $E$ is $p$-adically closed or
has a Henselian discrete valuation with a quasifinite residue field.
Note further that the statement of Theorem \ref{theo1.2} coincides in case
(1.2) (ii) with the local reciprocity law of the Artin-Tate abstract
class field theory \cite{AT}, Ch. 14, Sect. 5. As to Theorem \ref{theo1.3}, it
can be viewed as a partial analogue to the local Hasse symbol
(compare with Theorem \ref{theo3.2}, Remark \ref{rema7.1} and 
\cite{Wh1}, Theorem~3). The question of whether fields $E$ with LCFT 
are strictly PQL is open. By Theorem \ref{theo1.1}, its answer depends on the 
solution to one of the leading unsolved problems in Brauer group 
theory (see Remark \ref{rema3.4}, \cite{MS}, Sect. 16, and \cite{M}, 
Sect. 5). This allows us to prove in Section 3 that SQL-fields are 
those whose finite extensions admit LCFT. Note also that when 
$_{p}{\rm Br}(E)$ is finite, for each $p \in P(E)$, the answer is 
positive if and only if $E ^{\ast }/N(M/E)$ $\cong {\mathcal G}(M/E) 
^{{\rm Br}(E)}$, $M \in \Omega (E)$ (apply Theorem \ref{theo1.2} and 
\cite{P}, Sect. 15.1, Proposition~b).

\medskip
The paper is organized as follows: Section 2 includes generalities
about PQL-fields. The former and the latter assertions of Theorem
\ref{theo1.1} are proved in Sections 4 and 3, respectively. Our main 
results on local $p$-CFT are stated in Section 3. Theorems 
\ref{theo1.2}, \ref{theo1.3} and these results are proved in Sections 
4, 6 and 7. Section 5 contains an interpretation of a part of local 
$p$-CFT in terms of Galois cohomology. It generalizes \cite{LLMS}, 
Theorem~1, as well as known relations between local fields and 
Demushkin groups (cf. \cite{S2}, Ch. II, Theorem~4), and the 
sufficiency part of the main results of \cite{MW1,MW2}.
\medskip

\section{Preliminaries on PQL-fields}

\medskip
The present research is based on the possibility to reduce the
study of norm groups of finite abelian extensions to the
special case of $p$-extensions. The reduction is obtained by applying
the following lemma (which can be deduced from Galois theory (see
e.g., \cite{L}, Ch. VIII) and \cite{Ch6}, II, Lemma~2.2).

\medskip
\begin{lemm}
\label{lemm2.1}
Let $E$ be a field, $M \in \Omega (E)$, $M \neq E$, $\Pi $ the set
of prime divisors of $[M\colon E]$, and $M _{p} = M \cap E (p)$,
for each $p \in \Pi $. Then $M$ coincides with the compositum of the
fields $M _{p}\colon \ p \in \Pi $, $N(M/E) = \cap _{p \in \Pi } N(M
_{p}/E)$ and $E ^{\ast }/N(M/E)$ is isomorphic to the direct product
$\prod _{p \in \Pi } E ^{\ast }/N(M _{p}/E)$.
\end{lemm}

\medskip
The main results of \cite{Ch6}, I, used in this paper can be
stated as follows:

\medskip
\begin{prop}
\label{prop2.2}
Let $E$ be a $p$-quasilocal field with {\rm Br}$(E) _{p}
\neq \{0\}$, for some $p \in P(E)$. Suppose further that $R$ is a
finite extension of $E$ in $E (p)$ and $D \in d(E)$ is an algebra of
$p$-primary dimension. Then:
\par
{\rm (i)} $R$ is $p$-quasilocal and $D/E$ is a cyclic of exponent
{\rm exp}$(D) = {\rm ind}(D)$;
\par
{\rm (ii)} Br$(R) _{p}$ is a divisible group unless $p = 2$, $R = E$
and $E$ is formally real; when $E$ is formally real, $E (2) =
E(\sqrt{-1})$ and {\rm Br}$(E) _{2}$ is of order $2$;
\par
{\rm (iii)} $\rho _{E/R}$ maps {\rm Br}$(E) _{p}$ surjectively on
{\rm Br}$(R) _{p}$ and {\rm Cor}$_{R/E}$ maps {\rm Br}$(R) _{p}$
injectively in {\rm Br}$(E) _{p}$; every $E$-automorphism $\psi $ of
the field $R$ is extendable to a ring automorphism on each $D _{R}
\in s(R)$ of $p$-primary index;
\par
{\rm (iv)} $R$ embeds in $D$ as an $E$-subalgebra if and only
if the degree $[R\colon E]$ divides {\rm ind}$(D)$; $R$ is a
splitting field of $D$ if and only if ${\rm ind}(D) \mid [R\colon E]$.
\end{prop}

\medskip
Our next result gives an equivalent form of Proposition \ref{prop2.2} (iv),
for PQL-fields, and shows its optimality in the class of algebraic
strictly PQL-extensions of the field $\mathbb Q$ of rational
numbers. We refer the reader to \cite{Ch6}, II, for a proof of this
result, which demonstrates the applicability of the arithmetic
method of constructing such extensions, based on \cite{Ch9},
Theorem~2.2.

\medskip
\begin{prop}
\label{prop2.3}
{\rm (i)} Let $E$ be a {\rm PQL}-field, $M/E$ a finite
Galois extension and $R$ an intermediate field of $M/E$. If ${\mathcal
G}(M/E)$ is nilpotent, then $R$ embeds as an $E$-subalgebra in each
$D \in d(E)$ of index divisible by $[R\colon E]$.
\par
{\rm (ii)} For each nonnilpotent finite group $G$, there exists a
strictly {\rm PQL}-field $\Psi = \Psi (G)$, such that $E/\mathbb Q$
is an algebraic extension, {\rm Br}$(\Psi ) \cong \mathbb Q/\mathbb
Z$ and there is a Galois extension $\Psi ^{\prime }$ of $\Psi $ with
$\mathcal{G}(\Psi ^{\prime }/\Psi ) \cong G$, which does not embed
as a $\Psi $-subalgebra in any $\Delta \in d(\Psi )$ of index {\rm
ind}$(\Delta ) = [\Psi ^{\prime }\colon \Psi ]$.
\end{prop}

\medskip
The main results of \cite{Ch6}, I, and \cite{Ch2}, Sect. 3, show that
Br$(E)$ and $\mathcal{G} _{E}$ are the main algebraic structures
associated with any quasilocal field $E$. Therefore, we would like
to point out that $\mathcal{G} _{E}$ is prosolvable and Br$(E)$ is
embeddable in $\mathbb Q/\mathbb Z$ in the following two cases:

\medskip
(2.1) (i) $E$ is an algebraic extension of a quasilocal field $K$
with a Henselian valuation $v$, such that $v(K)$ is totally
indivisible (i.e. $v(K) \neq pv(K)$, for every $p \in \mathbb P$);
then Br$(E)$ is divisible with Br$(E) _{p'} = \{0\}$, $p ^{\prime }
\notin P(E)$ (concerning $\mathcal{G} _{E}$ and Br$(E)$, see
\cite{Ch1}, I, Lemma~1.2 and Proposition~3.1, and \cite{Ch11},
Corollary~5.3 as well as \cite{Ch6}, I, (1.3), respectively).
\par
(ii) $E$ is formally real and quasilocal; then Br$(E)$ is of order
$2$ and $\mathcal{G} _{E}$ has an abelian open subgroup of index
$2$, which is either procyclic or $2$-generated as a topological group
(see \cite{Ch2}, Sect. 3, and \cite{Ch10}, Proposition~3.4).

\medskip
\noindent Condition (2.1) (i) holds, for example, when $(K, v)$ is
Henselian discrete valued quasilocal or $E$ is an algebraic
quasilocal nonreal extension of a global field, such that Br$(E)
\neq \{0\}$ (cf. \cite{Ch9}, Sect. 3). The existence of a PQL-field
$F$ of essentially nonarithmetic nature, i.e. with Br$(F)$ not
embeddable in $\mathbb Q/\mathbb Z$, has been established in
\cite{SVdB}, Sect. 3, by observing that every abelian torsion group
$A$ embeds in Br$(F(A))$, for a suitably chosen strictly PQL-field
$F(A)$. This result and Proposition \ref{prop2.2} (ii) are complemented by the
following statement which describes, in conjunction with \cite{F},
Theorem~23.1, the isomorphism classes of Brauer groups of nonreal
PQL-fields (see \cite{Ch10}, Theorem~1.2 (i)-(ii), and for the
formally real case, \cite{Ch10}, Proposition~6.4). This statement
also shows that the absolute Galois groups of nonreal SQL-fields
need not be prosolvable and may have a very complex structure:

\medskip
(2.2) An abelian torsion group $T$ is isomorphic to Br$(E)$, for
some nonreal PQL-field $E = E(T)$ if and only if $T$ is divisible.
If $T$ is divisible, $S$ is a set of finite groups, $E _{0}$ is an
arbitrary field, and $T _{0}$ is a subgroup of Br$(E _{0})$
embeddable in $T$, then $E$ can be chosen so as to satisfy
the following:
\par
(i) $E$ is quasilocal, $P(E) = \mathbb P$ and $\rho _{E/L}$
surjective, for any finite extension $L/E$ (apply also the
Albert-Hochschild theorem (cf. \cite{S2}, Ch. II, 2.2)).
\par
(ii) $E/E _{0}$ is an extension, such that $E _{0}$ is algebraically
closed in $E$, $T _{0} \cap {\rm Br}(E/E _{0}) = \{0\}$ and every $G
\in S$ is realizable as a Galois group over $E$.

\medskip
When the set $_{p}T = \{t \in T\colon \ pt = 0\}$ is infinite, for
each $p \in \mathbb P$, (2.2) implies the following result (combined
with \cite{Ch11}, Remark~6.6, it points up the need for the
restrictions on the ground fields considered in \cite{KraMcK} and
\cite{RR}):

\medskip
(2.3) In order that $A, B \in d(E)$ have a common set of splitting
fields among the intermediate fields of $E _{\rm sep}/E$ it is
necessary and sufficient that ind$(A) = {\rm ind}(B)$ \cite{Ch6}, I,
Corollary~8.5. For each $n \in \mathbb N$, $d(E)$ contains
infinitely many nonisomorphic $E$-algebras of index $n$.

\medskip
Statements (2.1) and (2.2) can be supplemented by the following
partial conversion of the Koenigsmann-Neukirch theorem \cite{Koe},
Theorem~B:

\medskip
(2.4) If $(E, v)$ is a Henselian quasilocal field, such that all
finite groups are isomorphic to subquotients of $\mathcal{G} _{E}$
(by open normal subgroups), then $v(E)$ is divisible and every $D
\in d(E)$ is inertial relative to $v$ \cite{Ch11}, Proposition~6.3.

\medskip
\begin{rema}
\label{rema2.4}
Let $T$ be a divisible abelian torsion group and Supp$(T)$ the set
of those $p \in \mathbb P$ for which the $p$-component $T _{p}$ of
$T$ is nontrivial. Fix a field $E$ with Br$(E) \cong T$ as in (2.2)
(i) and denote by $E _{\rm sol}$ the maximal Galois extension of $E$
in $E _{\rm sep}$ with a prosolvable Galois group. Then $E _{\rm
sol}/E$ possesses an intermediate field $E ^{\prime }$ that is
strictly PQL with Br$(E ^{\prime }) \cong T$ (take as $E ^{\prime
}$, e.g., the fixed field of some Hall pro-${\rm Supp}(T)$-subgroup
of $\mathcal{G}(E _{\rm sol}/E)$). Note also that if Supp$(T) =
\mathbb P$, then $E$ is SQL.
\end{rema}

\medskip
Remark \ref{rema2.4} and statements (1.2) (i) and (2.1)$\div $(2.4) draw
interest in the following open questions:

\medskip
(2.5) (i) Describe the isomorphism classes of Brauer groups of
SQL-fields.
\par
(ii) Let $E$ be an SQL-field with $\mathcal{G} _{E}$ prosolvable. 
Find whether the nonzero groups Br$(E) _{p}$ except, possibly, one 
are isomorphic to the quasicyclic $p$-group $\mathbb Z (p ^{\infty 
})$. Prove whether the Sylow pro-$p$-subgroups of $\mathcal{G} _{E}$ 
are of rank $r _{p} \le 2$, for all $p \in \mathbb P$, with at most 
$2$ exceptions.
\par
(iii) Let $(F, v)$ be a Henselian quasilocal field with Br$(F)
_{\pi }$ and Br$(E) _{\pi '}$ not embeddable in $\mathbb Q/\mathbb
Z$, for some $\pi , \pi ^{\prime } \in \mathbb P$. Is $v(F)$
divisible?

\medskip\noindent
The questions in (2.5) (ii) are related by Galois cohomology (cf.
\cite{Ta}, page 265, \cite{MS}, (11.5), \cite{S2}, Ch. I, 4.2, and
\cite{W2}, Lemma~7) and the results of Section 5. Also, by
\cite{Ch11}, (1.4), the assumptions of (2.1) (i) require that $r
_{p} \le 2$ whenever $p \neq {\rm char}(\widehat E)$. In view of
\cite{Koe}, Theorem~B, one can prove that a negative answer to the
latter question in (2.5) (ii) would imply the existence of a field
$F$ with Br$(F) \neq \{0\}$ whose nontrivial valuations have
Henselizations equal to $F _{\rm sep}$. Moreover, if $3 \le r _{p}
\le r _{p'} \in \mathbb N$, for some $p, p ^{\prime } \in \mathbb
P$, $p < p ^{\prime }$, this would solve affirmatively Problem
11.5.9 (b) in \cite{FJ}.

\medskip
\begin{rema}
\label{rema2.5}
It is known (see \cite{NP}) that if $E$ is a field, then the triple
AT$(E)$ defined in (1.2) (ii) is an Artin-Tate class formation if
and only if $E$ is SQL, Br$(E)$ is embeddable in $\mathbb Q/\mathbb
Z$, and $\rho _{E/R}$ is surjective, for every finite extension $R$
of $E$ in $E _{\rm sep}$. Observing that (2.2) provides access to
the richest presently known sources of such fields, we add to the
examples given after the statement of (1.2) that AT$(E)$ is
Artin-Tate when $E$ is an algebraic SQL-extension of a global field
(see \cite{NP} and \cite{Ch9}, Sect. 3).
\end{rema}

\medskip

\section{Statements of the main results of local $p$-class field theory}

\medskip
Our main results on local $p$-CFT are stated as the following
two theorems:

\medskip
\begin{theo}
\label{theo3.1}
Let $E$ be a $p$-quasilocal field with {\rm Br}$(E)
_{p} \neq \{0\}$, for some $p \in P(E)$, and let $\Omega _{p} (E)$
and $d(p)$ be as in the Introduction. Then $E$ admits local $p$-{\rm
CFT} and, for each $M \in \Omega _{p} (E)$, $E ^{\ast }/N(M/E)$ is
isomorphic to the group $\mathcal{G}(M/E) ^{d(p)}$ defined in
Theorem \ref{theo1.2}.
\end{theo}

\medskip
Theorem \ref{theo3.1} plays a major role in the proof of the 
embeddability of Br$(E)$ into $\mathbb Q/\mathbb Z$ in case (2.1) 
(i), which in turn enables one to characterize the quasilocal 
property in the class of Henselian valued fields with totally 
indivisible value groups (see \cite{Ch11}, Sect. 6, and \cite{Ch4}, 
II). Note also that Theorem \ref{theo3.1} and Lemma \ref{lemm2.1} imply 
Theorem \ref{theo1.2} and the former assertion of Theorem 
\ref{theo1.1}. In addition, Theorem \ref{theo3.1} contains exact 
analogues to the fundamental correspondence and the local 
reciprocity law. Similarly, our next result can be viewed as such an 
analogue to Hasse's symbol for local fields (see \cite{Ko2}, Ch. 2, 
as well as Remark \ref{rema7.1} and Corollary \ref{coro7.3} below).

\medskip
\begin{theo}
\label{theo3.2} 
Under the hypotheses of Theorem \ref{theo3.1}, let $E _{\infty }$
be a $\mathbb Z _{p}$-extension of $E$ in $E (p)$, and for any $n
\in \mathbb N$, let $\Gamma _{n}$ be the intermediate field of $E
_{\infty }/E$ of degree $[\Gamma _{n}\colon E] = p ^{n}$. Then there
exist sets
\par\noindent
$H _{p} (E ^{\prime }) = \{( \ , M ^{\prime }/E ^{\prime
})\colon \ E ^{\prime \ast } \to \mathcal{G}(M ^{\prime }/E ^{\prime
}) ^{d(p)}, M ^{\prime } \in \Omega _{p} (E ^{\prime })\}$,  $E
^{\prime } \in \Omega _{p} (E)$,
\par\noindent 
of surjective group homomorphisms satisfying the following:
\par
{\rm (i)} The kernel of $( \ , M ^{\prime }/E ^{\prime })$ is equal
to $N(M ^{\prime }/E ^{\prime })$, for each $E ^{\prime } \in
\Omega _{p} (E)$, $M ^{\prime } \in \Omega _{p} (E ^{\prime })$;
\par
{\rm (ii)} If $M \in \Omega _{p} (E)$ and $K$ is an intermediate
field of $M/E$, then $( \ , K/E) = \pi _{M/K} \circ ( \ , M/E)$, 
where $\pi _{M/K}\colon \ \mathcal{G}(M/E) ^{d(p)} \to 
\mathcal{G}(K/E) ^{d(p)}$ is the homomorphism acting componentwise 
as the natural mapping of $\mathcal{G}(M/E)$ on $\mathcal{G}(K/E)$;
\par
{\rm (iii)} In the setting of (ii), $(\lambda , M/K) = (N _{E} ^{K} 
(\lambda ), M/E)$, for each $\lambda \in K ^{\ast }$;
\par
{\rm (iv)} The maps $( \ , \Gamma _{n}/E)$, $n \in \mathbb N$, 
determine the sets $H _{p} (E ^{\prime }), E ^{\prime } \in \Omega 
_{p} (E)$.
\end{theo}

\medskip
The rest of this Section concerns the open question of whether
fields with LCFT are strictly PQL. Lemma  \ref{lemm2.1} and first result in
this direction imply the latter assertion of Theorem \ref{theo1.1}.

\medskip
\begin{prop}
\label{prop3.3}
Let $E$ be a field admitting local $p$-{\rm CFT}, for some $p \in
P(E)$, and let $L$ be a degree $p$ extension of $E$ in $E (p)$. Then
{\rm Br}$(L/E) \neq \{0\}$ and {\rm Br}$(L/E)$ does not depend on
the choice of $L$. Furthermore, if each $D \in d(E)$ of exponent $p$
is similar to tensor products of cyclic division $E$-algebras of
index $p$, then $E$ is $p$-quasilocal.
\end{prop}

\begin{proof}
Our assumptions ensure that $L/E$ is cyclic, whence
Br$(L/E) \cong E ^{\ast }/N(L/E)$ (cf. \cite{L}, Ch. I, Sect. 6, and
\cite{P}, Sect. 15.1, Proposition~b). As $E$ admits local $p$-CFT
and $L \neq E$, this yields $N(L/E) \neq N(E/E) = E ^{\ast }$ and
Br$(L/E) \neq \{0\}$. In order to complete our proof, it suffices to
show that $L$ embeds as an $E$-subalgebra in each cyclic division
$E$-algebra of index $p$. If $\mathcal{G}(E (p)/E)$ is procyclic, this
is evident, so we assume that $\Omega _{p} (E)$ contains a field $F
\neq L$, such that $[F\colon E] = p$. It follows from Galois theory that
$LF \in \Omega (E)$, $[LF\colon E] = p ^{2}$ and $\mathcal{G}(LF/E)$ is
noncyclic. Therefore, $LF/E$ possesses $p + 1$ intermediate fields
of degree $p$ over $E$. Let $F ^{\prime }$ be such a field different
from $L$ and $F$. Then $N(F/E)N(F ^{\prime }/E) = E ^{\ast }$ and
$LF ^{\prime } = LF$. Moreover, the norm maps $N _{E} ^{F}$ and
$N _{E} ^{F'}$ are induced by $N _{L} ^{LF}$, so it turns out that $E
^{\ast } \subseteq N(LF/L)$. Hence, by \cite{P}, Sect. 15.1,
Proposition~b, $L$ embeds over $E$ in each $\Delta \in d(E)$ split
by $F$, which proves Proposition \ref{prop3.3}.
\end{proof}

\medskip
\begin{rema}
\label{rema3.4}
Let $E$ be a field and $p \in P(E)$.
\par
(i) If $\mathcal{G}(E(p)/E) \cong \mathbb Z _{p}$ and $L$ is the
unique degree $p$ extension of $E$ in $E(p)$, then $E$ admits local
$p$-CFT if and only if Br$(L/E) \neq \{0\}$ (apply Proposition \ref{prop3.3}
and \cite{P}, Sect. 15.1, Corollary~b).
\par
(ii) When Br$(E) _{p} \neq \{0\}$, the concluding condition of
Proposition \ref{prop3.3} is satisfied in the following cases: (a) $E$ is an
algebraic extension of a global or local field $E _{0}$; (b) $E$
contains a primitive $p$-th root of unity or char$(E) = p$; (c) $p =
3, 5$. In cases (b) and (c), this follows from the Merkur'ev-Suslin
theorem \cite{MS}, (16.1), \cite{A}, Ch. VII, Theorem~28, \cite{M},
Sect. 4, Corollary, and \cite{Ma}. In case (a), by class field
theory, every $A _{0} \in s(E _{0})$ is cyclic (cf. \cite{AT}, Ch.
10, Theorem~5), which implies the same, for all $A \in s(E)$.
\par
(iii) It is not known whether $E$ is strictly PQL, if it admits LCFT
and has a Henselian discrete valuation (see \cite{Ch4}, Sect. 2).
\end{rema}

\medskip
\begin{coro}
\label{coro3.5}
A finite purely inseparable extension $K$ of a field $E$ of
characteristic $p > 0$ admits local $p$-{\rm CFT} if and only if so
does $E$.
\end{coro}

\begin{proof}
By \cite{Ch6}, I, Proposition~4.4, $E$ is $p$-quasilocal if and only 
if $K$ is $p$-quasilocal, and by the Albert-Hochschild theorem, 
$\rho _{E/K}$ is surjective. Since the exponent of Br$(K/E)$ divides 
$[K\colon E]$ (see \cite{P}, Sects. 13.4 and 14.4), and by Witt's 
theorem (cf. \cite{Dr}, page 110), Br$(E) _{p}$ and Br$(K) _{p}$ are 
divisible, this implies that Br$(E) _{p} \neq \{0\} \leftrightarrow 
{\rm Br}(K) _{p} \neq \{0\}$, so Corollary \ref{coro3.5} follows 
from Theorem \ref{theo3.1} and Remark \ref{rema3.4} (ii).
\end{proof}

\medskip
The concluding result of this Section clarifies the role of
SQL-fields in LCFT. In view of (2.2), Remark \ref{rema2.5} and the main
results of \cite{Ch10}, it also determines the place in the study of
strictly PQL-fields of the Neukirch-Perlis variant of the theory
\cite{NP}, built upon (1.2) (ii).

\medskip
\begin{prop}
\label{prop3.6}
A field $E$ is {\rm SQL} if and only if its finite
extensions admit {\rm LCFT}. When occurs, {\rm Br}$(E _{p}) \neq
\{0\}$, provided that $E _{p}$ is the fixed field of a Sylow
pro-$p$-subgroup of $\mathcal{G}_{E}$, where $p \in \mathbb P$ is
chosen so that $\mathcal{G}_{E}$ is of cohomological $p$-dimension
{\rm cd}$_{p} (\mathcal{G} _{E}) \neq 0$.
\end{prop}

\begin{proof}
The left-to-right implication is contained in Theorem \ref{theo1.1}, 
so we prove the converse one and the nontriviality of Br$(E _{p}) 
_{p}$, for any $p \in \mathbb P$ satisfying cd$_{p} (\mathcal{G} 
_{E}) \neq 0$. Suppose that every finite extension $L$ of $E$ admits 
LCFT, and $B _{p}$ is the extension of $E$ generated by the $p$-th 
roots of unity in $E _{\rm sep}$. It is known (cf. \cite{L}, Ch. 
VIII, Sect. 3) that $B _{p} \subseteq E _{p}$, i.e. $E _{p}$ 
contains a primitive $p$-th root of unity unless $p = {\rm 
char}(E)$. Also, Proposition \ref{prop3.3} and Remark \ref{rema3.4} 
(ii) indicate that $L$ is $p$-quasilocal, provided that $B _{p} 
\subseteq L$. As $p \dagger [L ^{\prime }\colon E]$, for any finite 
extension $L ^{\prime }$ of $E$ in $E _{p}$, this enables one to 
obtain from general properties of scalar extension maps and Schur 
indices (cf. \cite{P}, Sect. 13.4, and \cite{Ch6}, I, (1.3)) that $E 
_{p}$ is $p$-quasilocal. Moreover, it follows from \cite{Ch6}, I, 
Lemma~8.3, and the choice of $p$ that $E$ is quasilocal. Since, by 
Proposition 3.3, PQL-fields with LCFT are strictly PQL, the obtained 
result shows that $E$ is SQL. Note further that the inequality 
cd$_{p}(\mathcal{G} _{E}) \neq 0$ ensures the existence of finite 
Galois extensions of $E$ in $E _{\rm sep}$ of degrees divisible by 
$p$. Therefore, by Sylow's theorems and Galois theory, there is a 
finite extension $ L _{0}$ of $E$ in $E _{p}$, such that $p \in P(L 
_{0})$. As in the proof of the statement that $E _{p}$ is 
$p$-quasilocal, these observations imply $p \in P(E _{p})$ and Br$(E 
_{p}) \neq \{0\}$.
\end{proof}

\medskip
Let us note that the abstract approach to LCFT dates back to the 
early 1950's (cf. \cite{AT}, Sect. 14). It fits the character of the 
related class formation theory, see \cite{K}, and accounts for the 
fact that the Neukirch-Perlis variant of LCFT goes substantially 
beyond the limits of (2.1). The fact itself is established by 
summing up Proposition \ref{prop3.6}, Remark \ref{rema2.5}, 
statements (2.1), (2.2) and (2.4), and the results on (2.5) (ii) 
mentioned in Section 2.
\medskip

\section{Proof of Theorem \ref{theo3.1}}

\medskip
Let $E$ be a $p$-quasilocal field, for a given $p \in P(E)$. By
Proposition \ref{prop2.2} (iv), then Br$(M/E) = \{b \in {\rm 
Br(E)}\colon \ [M\colon E]b = 0\}$, for every cyclic $p$-extension 
$M/E$. Hence, by the structure of divisible abelian torsion groups, 
and by the fact that Br$(M/E) \cong E ^{\ast }/N(M/E)$ (cf. 
\cite{F}, Theorem~23.1, and \cite{P}, Sect. 15.1, Proposition~b), $E 
^{\ast }/N(M/E) \cong {\mathcal G}(M/E) ^{d(p)}$. To prove the 
obtained isomorphism, for any $M \in \Omega _{p} (E)$, we need the 
following lemmas.

\medskip
\begin{lemm}
\label{lemm4.1}
Let $E$ be a $p$-quasilocal field, for some $p \in P(E)$, and let $E 
_{1}, \dots E _{t}$ be cyclic extensions of $E$ in $E (p)$, for a 
given integer $t \geq 2$. Assume that the compositum $E ^{\prime }$ 
of the fields $E _{j}\colon \ j = 1, \dots , t$, satisfies the 
equality $[E ^{\prime }\colon E] = \prod _{j=1} ^{t} [E _{j}\colon 
E]$. Then $N(E ^{\prime }/E) = \bigcap _{j=1} ^{t} N(E _{j}/E)$.
\end{lemm}

\begin{proof}
The inclusion $N(E ^{\prime }/E) \subseteq \bigcap _{j=1} ^{t} N(E 
_{j}/E)$ follows from the transitivity of norm maps in towers of
finite separable extensions (cf. \cite{L}, Ch. VIII, Sect. 5).
Conversely, let $c \in \bigcap _{j=1} ^{t} N(E _{j}/E)$ and $\b \in
E _{1} ^{\ast }$ be of norm $N _{E} ^{E _{1}} (\b ) = c$. The 
equality $[E ^{\prime }\colon E] = \prod _{j=1} ^{t} [E _{j}\colon 
E]$, Galois theory and \cite{Ch3}, Lemma~4.2, imply $[E ^{\prime 
}\colon E _{1}]  = \prod _{i=2} ^{t} [(E _{1}E _{i})\colon E _{1}]$ 
and $\b \in \bigcap _{i=2} ^{t} N(E _{1}E _{i}/E _{1})$. This proves 
Lemma \ref{lemm4.1} in the case where $t = 2$. Since $(E _{1}E _{i})/E _{1}$ 
is cyclic, for each $i \geq 2$, and by Proposition \ref{prop2.2} 
(i), $E _{1}$ is $p$-quasilocal, the obtained result makes it easy 
to complete our proof by induction on $t$.
\end{proof}

\medskip
\begin{lemm}
\label{lemm4.2}
Assume that $E$, $F$ and $L$ are fields, such that $E$ is 
$p$-quasilocal, $L \in \Omega _{p} (E)$, $E \subseteq F \subseteq L$ 
and $F/E$ is cyclic. Then $\psi (\a)\alpha ^{-1} \in N(L/F)$, for 
each $\alpha \in F ^{\ast }$ and $\psi \in \mathcal{G}(F/E)$.
\end{lemm}

\begin{proof}
As $F$ is $p$-quasilocal and $L \in \Omega _{p} (F)$, whence 
${\mathcal G}(L/F)$ decomposes into a direct product of cyclic 
groups, Galois theory and Lemma \ref{lemm4.1} allow one to consider only the
case in which $L/F$ is cyclic. Let $\psi ^{\prime }$ be an 
automorphism of $L$ extending $\psi $. Fix a generator $\sigma $ of 
$\mathcal{G}(L/F)$, denote by $A _{\alpha }$ the cyclic $F$-algebra 
$(L/F, \sigma , \alpha )$, for an arbitrary $\alpha \in F ^{\ast }$, 
put $m = [L\colon F]$, and take an invertible element $\eta \in A 
_{\alpha }$ so that $\eta ^{m} = \alpha $ and $\eta \lambda \eta 
^{-1} = \sigma (\lambda )$, for every $\lambda \in L$. Then 
Proposition \ref{prop2.2} (iii) and the Skolem-Noether theorem (cf. 
\cite{P}, Sect. 12.6) imply that $A _{\alpha }$ has a ring 
automorphism $\tilde \psi $, such that $\tilde \psi (\lambda )$ $= 
\psi ^{\prime } (\lambda )$, for any $\lambda \in L$. Since $\sigma 
\psi ^{\prime } = \psi ^{\prime }\sigma $, this ensures that the 
element $\eta ^{-1}\tilde \psi (\eta ) := \mu $ lies in the 
centralizer of $L$ in $A _{\alpha }$. Thereby, we have $\mu \in L 
^{\ast }$ and $\tilde \psi (\eta ) ^{m} = \psi (\alpha ) = \eta 
^{m}N _{F} ^{L} (\mu ) = \alpha N _{F} ^{L} (\mu )$, which proves 
Lemma \ref{lemm4.2}.
\end{proof}

\medskip
Assuming that $E$ is a $p$-quasilocal field, $M _{u} \in \Omega
_{p} (E)$, $[M _{u}\colon E] = p ^{\mu _{u}}\colon \ u = 1, 2$, and
putting $M ^{\prime } =  M _{1}M _{2}$, we prove the following
assertions:

\medskip
(4.1) (i) If $N(M _{1}/E) = N(M _{2}/E)$ and $M _{1} \subseteq M
_{2}$, then $M _{1} = M _{2}$;
\par
(ii) If $M _{1} \cap M _{2} = E$, then $N(M _{1}/E) \cap N(M _{2}/E)
= N(M ^{\prime }/E)$,
\par
\noindent $N(M _{1}/E)N(M _{2}/E) = E ^{\ast }$, each $\sigma \in
\mathcal{G}(M _{2}/E)$ has a unique prolongation $\sigma ^{\prime } \in
\mathcal{G}(M ^{\prime }/M _{1})$, and the mapping of $\mathcal{G}(M
_{2}/E)$ on $\mathcal{G}(M ^{\prime }/M _{1})$ by the rule $\sigma \to
\sigma ^{\prime }$ is an isomorphism;
\par
(iii) Under the hypotheses of (ii), if $M _{2}/E$ is cyclic and
$\mathcal{G}(M _{2}/E) = \langle \sigma \rangle $, then $\mathcal{G}(M
^{\prime }/M _{1}) = \langle \sigma ^{\prime } \rangle $ and Cor$_{M
_{1}/E}$ maps Br$(M ^{\prime }/M _{1})$ into Br$(M _{2}/E)$ by the
formula $[(M ^{\prime }/M _{1}, \sigma ^{\prime }, \theta )] \to [(M
_{2}/E, \sigma , N _{E} ^{M _{1}} (\theta ))]$, $\theta \in M _{1} ^{\ast }$;
when Br$(E) _{p}$ is divisible, Cor$_{M _{1}/E}$ induces
isomorphisms Br$(M _{1}) _{p} \cong {\rm Br}(E) _{p}$ and Br$(M
^{\prime }/M _{1}) \cong {\rm Br}(M _{2}/E)$.

\medskip
It is sufficient to consider the special case where $M _{u} \neq
E\colon \ u = 1, 2$, and $M _{1} \neq M _{2}$. In view of
Proposition \ref{prop2.2} (ii), then Br$(E) _{p}$ is divisible. At the same
time, Galois theory and the assumptions on $M _{1}/E$ imply the
existence of a cyclic degree $p$ extension $M _{0}$ of $E$ in $M
_{1}$. Suppose that $M _{1} \subseteq M _{2}$ and $d \in N(M
_{u}/E)\colon \ u = 1, 2$, fix a generator $\psi $ of $\mathcal{G}(M
_{0}/E)$ and elements $\eta _{j} \in M _{j} ^{\ast }$, $j = 1, 2$,
so that $N _{E} ^{M _{1}} (\eta _{1}) = N _{E} ^{M _{2}} (\eta _{2})
= d$. It follows from Hilbert's Theorem 90 that $N _{M _{0}} ^{M
_{1}} (\eta _{1}) = N _{M _{0}} ^{M _{2}} (\eta _{2})\psi (\b )\b
^{-1}$, for some $\b \in M _{0} ^{\ast }$. Hence, by Lemma \ref{lemm4.2}, $N(M
_{1}/E) = N(M _{2}/E)$ if and only if $N(M _{1}/M _{0}) = N(M _{2}/M
_{0})$. Observing also that $M _{0}$ is $p$-quasilocal and $[M
_{1}\colon M _{0}] = p ^{\mu _{1}-1}$, and proceeding by induction
on $\mu _{1}$, one proves (4.1) (i).
\par

Assume now that $M _{1} \cap M _{2} = E$. Then the Galois-theoretic
parts of (4.1) (ii) and (iii) are contained in \cite{L}, Ch. VIII,
Theorem~4. The rest of the former assertion of (4.1) (iii) is
implied by Proposition \ref{prop2.2} (iii), the basic
restriction-corestriction (abbr, RC) formula (cf. \cite{T},
Theorem~2.5) and the lemmas in \cite{Ch3}, Sect. 4. The equality
$N(M ^{\prime }/E) = N(M _{1}/E) \cap N(M _{2}/E)$ follows from the
presentability of $M _{1}$ and $M _{2}$ as compositums of cyclic
extensions of $E$ satisfying the conditions of Lemma \ref{lemm4.1}. It remains
for us to show that $N(M _{1}/E)N(M _{2}/E) = E ^{\ast }$ and to
prove the latter part of (4.1) (iii). Suppose first that $M _{2}/E$
is cyclic and the automorphisms $\sigma , \sigma ^{\prime }$ are determined
as in (4.1) (iii), fix an element $c \in E ^{\ast }$ out of $N(M
_{2}/E)$, and put $A _{c} =  (M _{2}/E, \sigma , c)$. It is clear from
\cite{P}, Sect. 15.1, Proposition~b, and the choice of $c$ that
ind$(A _{c}) \mid p ^{\mu _{2}}$ and ind$(A _{c}) > 1$. As Br$(E)
_{p}$ is divisible, there exists $\Delta _{c} \in d(E)$, such that
$p ^{\mu _{1}}[\Delta _{c}] = [A _{c}]$. In addition, it follows
from Proposition \ref{prop2.2} (i) that ind$(\Delta _{c}) = p ^{\mu _{1}}.
{\rm ind}(A _{c})$. Observing also that $[M ^{\prime }\colon E] = p
^{\mu }$ and ind$(\Delta _{c}) \mid p ^{\mu }$, where $\mu = \mu
_{1} + \mu _{2}$, one obtains from Proposition \ref{prop2.2} (iv) that
$[\Delta _{c}] \in {\rm Br}(M ^{\prime }/E)$. Since $\Delta _{c}
\otimes _{E} M ^{\prime } \cong (\Delta _{c} \otimes _{E} M _{1})
\otimes _{M _{1}} M ^{\prime }$ as $M ^{\prime }$-algebras (cf.
\cite{P}, Sect. 9.4, Corollary~a), this means that $[\Delta
_{c} \otimes _{E} M _{1}] \in {\rm Br}(M ^{\prime }/M _{1})$, or
equivalently, that $[\Delta _{c} \otimes _{E} M _{1}] = [(M ^{\prime
}/M _{1}, \sigma ^{\prime }, \alpha )]$, for some $\alpha \in M _{1}
^{\ast }$. Therefore, by Proposition \ref{prop2.2} (iii), the RC-formula and
the former statement of (4.1) (iii), $[A _{c}] = [(M _{2}/E, \sigma
, N _{E} ^{M _{1}} (\alpha ))]$. As $[A _{c}\colon E] = [(M _{2}/E,
\sigma , N _{E} ^{M _{1}} (\alpha ))\colon E]$ $= [M _{2}\colon E]
^{2}$, this proves that $A _{c} \cong (M _{2}/E, \sigma , N _{E} ^{M
_{1}} (\alpha ))$ over $E$. Hence, by \cite{P}, Sect. 15.1,
Proposition~b, $c.N _{E} ^{M _{1}} (\alpha ) ^{-1} \in N(M _{2}/E)$,
which yields $N(M _{1}/E)N(M _{2}/E) = E  ^{\ast }$ in case $M
_{2}/E$ is cyclic. Our argument, combined with Proposition \ref{prop2.2}
(ii)-(iii) and the RC-formula, also proves the latter part of (4.1)
(iii). Henceforth, we assume that $M _{2}/E$ is noncyclic, i.e.
$\mathcal{G}(M _{2}/E)$ is an abelian $p$-group of rank $r \geq 2$.
Then it follows from Galois theory and the structure of finite
abelian groups that there exist cyclic extensions $F _{1}$ and $F
_{2}$ of $E$ in $M _{2}$ such that $F _{1} \cap F _{2} = E$ and the
$p$-groups $\mathcal{G}(M _{2}/F _{1})$ and $\mathcal{G}(M _{2}/F
_{2})$ are of rank $r - 1$. As $M _{1} \cap M _{2} = E$, one also
sees that $(M _{1}F _{u}) \cap M _{2} = F _{u}\colon \ u = 1, 2$.
Taking now into account that $F _{1}$ and $F _{2}$ are
$p$-quasilocal fields, and arguing by induction on $r$, one
concludes that it suffices to deduce the equality $N(M _{1}/E)N(M
_{2}/E) = E ^{\ast }$ under the extra hypothesis that $N(M _{1}F
_{u}/F _{u})N(M _{2}/F _{u}) = F _{u} ^{\ast }\colon \ u = 1, 2$.
Then, by norm transitivity in towers of finite extensions, $N(F
_{1}/E)N(F _{2}/E) \subseteq N(M _{1}/E)N(M _{2}/E)$, and since $N(F
_{1}/E)N(F _{2}/E) = E ^{\ast }$, this completes the proof of (4.1).

\medskip
\begin{rema}
\label{rema4.3}
Let $E$ be a field and $M _{1}$, $M _{2}$ be finite
extensions of $E$ in $E _{\rm sep}$, such that $M _{2}/E$ is cyclic
and $M _{1} \cap M _{2} = E$. It is known that then the former
statement of (4.1) (iii) remains valid. The assertions on ${\mathcal
G}(M _{2}/E)$ and $\mathcal{G}(M ^{\prime }/M _{1})$ follow from
\cite{L}, Ch. VIII, Theorem~4, and the formula for the action of
Cor$_{M _{1}/E}$ on Br$(M ^{\prime }/M _{1})$ can be proved by a
group-cohomological technique (cf. \cite{We}, Proposition~4.3.7). It
is therefore worth noting that (4.1) (iii) plays a role in the proof
not only of Theorems \ref{theo3.1} and \ref{theo3.2} but also of 
Proposition \ref{prop2.2} (i) (see \cite{Ch6}, I, Sect. 7). This, 
combined with \cite{Ch10}, Theorem~1.2 (i)-(ii), enables one to find 
alternative field-theoretic proofs of the formula in (4.1) (iii) and 
of other results on Cor$_{F/E}$, for any finite extension $F$ of $E$ 
in $E _{\rm sep}$ (by reduction to the setting of (2.2) (i), see 
\cite{Ch7}).
\end{rema}

\medskip
We are now in a position to prove Theorem \ref{theo3.1} (and 
thereby, Theorem \ref{theo1.2} and the former part of Theorem 
\ref{theo1.1} as well). Assuming as above that $M _{u} \in \Omega 
_{p} (E)\colon \ u = 1, 2$, and $M _{1}M _{2} = M  ^{\prime }$, put 
$L ^{\prime } = M _{1} \cap M _{2}$. We first show that $N(L 
^{\prime }/E) = N(M _{1}/E)N(M _{2}/E)$ and $N(M ^{\prime }/E) = N(M 
_{1}/E) \cap N(M _{2}/E)$. In view of (4.1) (ii), it suffices to 
consider the case where $L ^{\prime } \neq E$. As $L ^{\prime }$ is 
$p$-quasilocal, (4.1) (ii) yields $L ^{\prime \ast } = N(M _{1}/L 
^{\prime })N(M _{2}/L ^{\prime })$, so the equality $N(L ^{\prime 
}/E) = N(M _{1}/E)N(M _{2}/E)$ is obtained from norm transitivity in 
the towers $E \subset L ^{\prime } \subseteq M _{u}\colon \ u = 1, 
2$. Suppose further that $[L ^{\prime }\colon E] = p ^{m}$ and fix
a degree $p$ extension $L$ of $E$ in $L ^{\prime }$. As in the proof
of (4.1) (i), it is seen that an element $\lambda \in L ^{\ast }$ 
lies in $N(M _{i}/L)$, for some $i \in \{1, 2\}$, if and only if $N 
_{E} ^{L} (\lambda) \in N(M _{i}/E)$. Since $[L ^{\prime }\colon L] 
= p ^{m-1}$ and $L$ is $p$-quasilocal, this makes it easy to prove 
inductively that $N(M ^{\prime }/E) = N(M _{1}/E) \cap N(M _{2}/E)$. 
It follows from this result and (4.1) (i) that the natural mapping 
$\Omega _{p} (E) \to {\rm Nr}(E)$ is injective. Thus the statement 
that $E$ admits local $p$-CFT is proved, which allows us to deduce 
Theorem \ref{theo1.1} from Lemma \ref{lemm2.1} and Proposition \ref{prop3.3}.
\par
To finish the proof of Theorem \ref{theo3.1} (and Theorem 
\ref{theo1.2}) we show that $E ^{\ast }/N(M/E) \cong 
\mathcal{G}(M/E) ^{d(p)}$, provided that $M \in \Omega _{p} (E)$ and 
$\mathcal{G}(M/E)$ is noncyclic. Then $\mathcal{G}(M/E)$ is an 
abelian $p$-group of rank $r(M) \geq 2$, so it follows from Galois 
theory that $M/E$ has intermediate fields $\Phi _{1}$ and $\Phi 
_{2}$, such that $\Phi _{1}\Phi _{2} = M$, $\Phi _{1} \cap \Phi _{2} 
= E$ and $\Phi _{2}/E$ is cyclic. This means that $\mathcal{G}(M/E)$ 
is isomorphic to the direct product $\mathcal{G}(\Phi _{1}/E) \times 
\mathcal{G}(\Phi _{2}/E)$, and $\mathcal{G}(\Phi _{1}/E)$ is of rank 
$r(M) - 1$ as a $p$-group. Since $N(\Phi _{1}/E)N(\Phi _{2}/E) = E 
^{\ast }$ and $N(\Phi _{1}/E) \cap N(\Phi _{2}/E) = N(M/E)$, the 
natural diagonal embedding of $E ^{\ast }$ into $E ^{\ast } \times E 
^{\ast }$ induces a group isomorphism $E ^{\ast }/N(M/E) \cong E 
^{\ast }/N(\Phi _{1}/E) \times E ^{\ast }/N(\Phi _{2}/E)$. Now our 
proof is easily completed proceeding by induction on $r(M)$.

\medskip
\begin{rema}
\label{rema4.4}
Given a PQL-field $E$, denote by ${\rm c}(M)$ the intersection of
fields $N \in \Omega (E)$ with $N(N/E) = N(M/E)$, for each $M \in
\Omega (E)$. It follows from Theorem \ref{theo3.1}, Lemma \ref{lemm2.1} and 
\cite{Ch6}, I, Lemma~4.2 (ii), that the natural map $\nu $ of 
$\Omega (E)$ on the set $N _{\rm ab} (E) = \{N(M/E)\colon \ M \in 
\Omega (E)\}$ satisfies (1.1). Therefore, Br$(E) _{p} \neq \{0\}$ 
and $p \dagger [M\colon c(M)]$ whenever $M \in \Omega (E)$ and $p 
\in \mathbb P$ divides $[{\rm c}(M)\colon E]$. Note also that 
$N({\rm c}(M)/E) = N(M/E)$ and $N(M _{0}/E) \neq N(M/E)$ in case $M 
_{0} \in \Omega (E)$, $M _{0} \subseteq {\rm c}(M)$ and $M _{0} \neq 
{\rm c}(M)$. Since the set Cl$(E) = \{{\rm c}(\Lambda )\colon \ 
\Lambda \in \Omega (E)\}$ is closed under taking subextensions of 
$E$ and finite compositums, whence $\nu $ induces a bijection of 
Cl$(E)$ on $N _{\rm ab} (E)$, these observations allow us to view 
${\rm c}(M)$ as a class field of $N(M/E)$.
\end{rema}
\medskip

\section{Galois cohomological interpretation of Theorem \ref{theo3.1}}

\medskip
In this Section we consider some Galois cohomological aspects of the
problem of characterizing fields with LCFT. Let $P$ an infinite
pro-$p$-group, cd$(P)$ the cohomological dimension of $P$ and
$\mathbb F _{p}$ a field with $p$ elements, for some $p \in \mathbb
P$. We say that $P$ is a $p$-group of Demushkin type, if the
(continuous) cohomology group homomorphism $\varphi _{h}\colon \ H
^{1} (P, \mathbb F _{p}) \to H ^{2}(P, \mathbb F _{p})$ mapping each
$g \in H ^{1} (P, \mathbb F _{p})$ into the cup-product $h \cup g$
is surjective, for every $h \in H ^{1} (P, \mathbb F _{p}) \setminus
\{0\}$. We call a degree of $P$ the dimension of $H ^{2} (P, \mathbb
F _{p})$ as an $\mathbb F _{p}$-vector space. The defined groups and
local $p$-CFT are related as follows:

\medskip
\begin{prop}
\label{prop5.1}
Let $E$ be a nonreal field containing a primitive $p$-th root of
unity, for some $p \in P(E)$. Then the following conditions are
equivalent:
\par
{\rm (i)} $E$ admits local $p$-{\rm CFT};
\par
{\rm (ii)} $\mathcal{G}(E (p)/E)$ is a $p$-group of Demushkin type
of degree $d \geq 1$;
\par
{\rm (iii)} $E$ is $p$-quasilocal with {\rm Br}$(E) _{p} \neq
\{0\}$;
\par
{\rm (iv)} ${\rm cd}(\mathcal{G}(E (p)/E)) = 2$ and $_{p} {\rm Br}(E
^{\prime })$ is a trivial module over the integral group ring
$\mathbb Z[\mathcal{G}(E ^{\prime }/E)]$, for every degree $p$
extension $E ^{\prime }$ of $E$ in $E (p)$.
\end{prop}

\begin{proof}
The equivalence (i)$\leftrightarrow $(iii) is implied by Proposition
\ref{prop3.3}, Remark \ref{rema3.4} (ii) and Theorem \ref{theo3.1}. 
Note also that Br$(E) _{p} = \{0\}$ if and only if 
$\mathcal{G}(E (p)/E)$ is a free pro-$p$-group, i.e. a $p$-group of 
Demushkin type of degree zero (cf. \cite{W1}, Theorem~3.1, \cite{W2}, 
page 725, and \cite{S2}, Ch. I, 4.1 and 4.2). In particular, this 
holds when $\mathcal{G}(E(p)/E)$ is of rank $1$ as a pro-$p$-group, 
since then $\mathcal{G}(E(p)/E) \cong \mathbb Z _{p}$ (see, e.g., 
\cite{Ch6}, I, Remark~3.4 (ii)). These observations, \cite{Ch6}, I, 
Lemma~3.8, and the end of Proposition \ref{prop2.2} (ii) prove that 
(ii)$\leftrightarrow $(iii). It remains to be seen that 
(iii)$\leftrightarrow $(iv). As $E$ is nonreal and $p \in P(E)$, one 
obtains from Galois theory and \cite{Wh2}, Theorem~2, that 
$\mathcal{G}(E(p)/E)$ possesses a closed normal subgroup $H$, such 
that $\mathcal{G}(E (p)/E)/H \cong \mathbb Z _{p}$. Since Br$(E) _{p} 
\neq \{0\}$, this means that $H \neq \{1\}$, so it follows from 
\cite{Ch6}, I, Proposition~4.6 (ii), and Galois cohomology (see 
\cite{S2}, Ch. I, 4.2 and Proposition~15) that cd$(H) = 1$ and 
cd$(\mathcal{G}(E (p)/E)) = 2$. Hence, by Teichm\"uller's theorem 
(cf. \cite{Dr}, Ch. 9, Theorem~4) and Proposition \ref{prop2.2} 
(iii), (iii)$\to $(iv). To prove that (iv)$\to $(iii) we need the 
following results (see \cite{MS}, (11.5), and \cite{Ko1}, 
Proposition~3.26):

\medskip
(5.1) If $F$ is a field containing a primitive $p$-th root of unity,
and $M$ is a finite Galois extension of $F$ in $F(p)$, then there
exists an isomorphism $\kappa _{M}\colon H ^{2} (\mathcal{G}(F
(p)/M), \mathbb F _{p}) \cong $ $_{p} {\rm Br}(M)$ as $\mathbb
Z[\mathcal{G}(M/F)]$-modules, and the compositions Cor$_{M/F} \circ
\kappa _{M}$ and $\kappa _{F} \circ $ ${\rm cor}_{M/F}$ coincide,
where cor$_{M/F}$ is the corestriction map $H ^{2}(\mathcal{G}(F
(p)/M), \mathbb F _{p}) \to H ^{2} (\mathcal{G}(F (p)/F), \mathbb F
_{p})$. Hence, Cor$_{M/F}$ is surjective if and only if so is
cor$_{M/F}$.

\medskip\noindent
It suffices to show that Br$(E ^{\prime }/E) =$ $_{p} {\rm Br}(E)$,
for an arbitrary degree $p$ extension $E ^{\prime }$ of $E$ in
$E(p)$. The equality cd$(\mathcal{G}(E (p)/E)) = 2$ ensures that
Br$(E) _{p} \neq \{0\}$ and it follows from (5.1) and \cite{NSW},
Proposition~3.3.8, that Cor$_{E'/E}$ maps $_{p} {\rm Br}(E ^{\prime
})$ surjectively on $_{p} {\rm Br}(E)$. We show that $_{p} {\rm
Br}(E ^{\prime })$ includes the preimage $\Pi _{p}(E)$ of $_{p} {\rm
Br}(E)$ in Br$(E ^{\prime }) _{p}$ under Cor$_{E'/E}$. Our
assumptions guarantee (in conjunction with Galois theory and the
normality of maximal subgroups of finite $p$-groups, see \cite{L},
Ch. VIII and Ch. I, Sect. 6) that $E ^{\prime }/E$ is cyclic.
Observe also that Br$(E ^{\prime }) _{p}$ is a trivial $\ZZ
[\mathcal{G}(E ^{\prime }/E)]$-module. Indeed, let $\sigma $ be a
generator of $\mathcal{G}(E ^{\prime }/E)$, and let $b \in {\rm
Br}(E ^{\prime })$ be an element of order $p ^{k}$, for some $k \in
\mathbb N$. Proceeding by induction on $k$, one sees that it
suffices to prove that $\sigma (b) = b$, provided that $k \ge 2$ and
$\sigma (b) = b + b _{0}$, for some $b _{0} \in _{p} {\rm Br}(E
^{\prime })$. The assumption on $b _{0}$ indicates that Cor$_{E'/E}
(b _{0}) = 0$, so it follows from (5.1) and \cite{LLMS}, Corollary,
that $b _{0} = a - \sigma (a)$, for some $a \in $ $_{p} {\rm Br}(E
^{\prime })$. This yields $\sigma (b + a) = b + a$, which implies
$\sigma (b) = b$ and our assertion about Br$(E ^{\prime }) _{p}$
(because $_{p} {\rm Br}(E ^{\prime })$ is a trivial $\ZZ
[\mathcal{G}(E ^{\prime }/E)]$-module). By Teichm\"uller's theorem,
the established property of Br$(E ^{\prime }) _{p}$ is equivalent to
the inclusion of Br$(E ^{\prime }) _{p}$ in the image of $\rho
_{E/E'}$. Since Cor$_{E'/E}$ maps $_{p} {\rm Br}(E ^{\prime })$
surjectively upon $_{p} {\rm Br}(E)$, it is now easy to obtain from
the RC-formula that $\Pi _{p}(E) \subseteq $ $_{p} {\rm Br}(E
^{\prime })$. Moreover, it becomes clear that the implication
(iv)$\to $(iii) will follow, if we show that $_{p} {\rm Br}(E)$ is
included in the subgroup $p{\rm Br}(E) _{p} = \{p\beta \colon \
\beta \in {\rm Br}(E) _{p}\}$ of Br$(E) _{p}$. Denote by $B$ the
extension of $E$ in $E (p)$ obtained by adjunction of a primitive $p
^{2}$-th root of unity. It follows from \cite{MS}, (16.1), \cite{P},
Sect. 15.1, Corollary~b, and Kummer theory that $_{p} {\rm Br}(B)
\subseteq p{\rm Br}(B) _{p}$. Hence, the surjectivity of the
homomorphism $_{p} {\rm Br}(B) \to $ $_{p} {\rm Br}(E)$ induced by
Cor$_{B/E}$ implies $_{p} {\rm Br}(E) \subseteq p{\rm Br}(E) _{p}$.
Therefore, Br$(E ^{\prime }/E) =$ $_{p} {\rm Br}(E)$ and (iv)$\to
$(iii), so Proposition \ref{prop5.1} is proved.
\end{proof}

\medskip
\begin{coro}
\label{coro5.2}
Let $E$ be an {\rm SQL}-field, such that ${\rm cd}_{p} (\mathcal{G}
_{E}) \neq 0$, for a given $p \in \mathbb P$. Then the Sylow
pro-$p$-subgroups of $\mathcal{G} _{E}$ are $p$-groups of Demushkin
type of degree $d _{p} \geq 1$, unless $p = {\rm char}(E)$ or $E$ is
formally real and $p = 2$.
\end{coro}

\begin{proof}
One may consider only the special case of $E  _{\rm sep} = E (p)$ 
(see the proof of Proposition \ref{prop3.6}). Then Br$(E) _{p} \neq 
\{0\}$ and our conclusion follows from Propositions \ref{prop5.1}, 
\ref{prop2.2} (ii) and \cite{S2}, Ch. II, Proposition~3.
\end{proof}

\medskip
It is easily seen that the degrees of $p$-groups of Demushkin type 
are bounded by their ranks. The following conversion of this fact
can be deduced from \cite{MS}, (11.5), Proposition \ref{prop5.1}, 
\cite{Ch6}, I, Theorem~8.1, and the sufficiency part of (2.2) (by 
specifying the cardinalities of the fields $E _{0}$ and $E$ in (2.2) 
(ii), see \cite{Ch10}, Remark~5.4, for more details):

\medskip
(5.2) For any system $d \geq \aleph _{0}$ and $d _{p} \leq d\colon \
p \in \mathbb P$, of cardinal numbers, there is a field $E$
containing a primitive $p$-th root of unity, for every $p \in 
\mathbb P$, and such that $\mathcal{G}(E (p)/E)$ and the Sylow
pro-$p$-subgroups of $\mathcal{G} _{E}$ are of rank $d$, Demushkin
type and degree $d _{p}$.

\medskip
By a Demushkin pro-$p$-group, we mean a $p$-group of Demushkin type
of degree $1$. Demushkin pro-$p$-groups of rank $r(P) \leq \aleph
_{0}$ have been classified by Demushkin, Labute and Serre (cf.
\cite{Lab1, Lab2} and further references there). When $r(P) = \aleph
_{0}$, by \cite{MW1} and \cite{MW2}, $P$ has $s$-invariant zero if
and only if $P \cong {\mathcal G}(F (p)/F)$, for some field $F$.
Hence, by applying Lemma~3.5 of \cite{Ch4} (in a modified form
adjusted to the case singled out by \cite{MW2}, Proposition~3.1
(iii)), and arguing as in the proof of (5.2), one supplements it as
follows:

\medskip
(5.3) For each sequence $G _{p}\colon \ p \in \mathbb P$, of
Demushkin pro-$p$-groups of rank $\aleph _{0}$ and $s$-invariant
zero, there exists a field $E$ such that $\mathcal{G}(E (p)/E)$ and
the Sylow pro-$p$-subgroups of $\mathcal{G} _{E}$ are isomorphic to
$G _{p}$ when $p$ ranges over $\mathbb P$.

\medskip
Next we show that the inequality $d \geq \aleph _{0}$ in (5.2) is 
essential. This result is a special case of \cite{Ku}, Corollary~5 
(see also \cite{LLMS}, Sects. 1 and 2). For convenience of the 
reader, we present it here with a short proof based on Theorem 
\ref{theo3.1} and Proposition \ref{prop5.1}.

\medskip
\begin{coro}
\label{coro5.3}
Let $E$ be a field containing a primitive $p$-th root of unity, for
some $p \in P(E)$, and with $\mathcal{G}(E (p)/E)$ a $p$-group of
Demushkin type of finite rank $r(p)$. Then $\mathcal{G}(E (p)/E)$ is
a Demushkin group or a free pro-$p$-group.
\end{coro}

\begin{proof}
In view of \cite{W2}, Lemma~7, and Proposition \ref{prop5.1} (with 
its proof), one may consider only the case where $r(p) \geq 3$ and 
Br$(E) _{p} \neq \{0\}$. Take a field $E ^{\prime } \in \Omega (E)$ 
so that $\mathcal{G}(E ^{\prime }/E)$ has exponent $p$ and order $p 
^{r(p)-1}$. By Theorem \ref{theo3.1} and \cite{MS}, (11.5), then $E 
^{\ast p} \subseteq N(E ^{\prime }/E)$ and $N(E ^{\prime }/E)$ is a 
subgroup of $E ^{\ast }$ of index $p ^{(r(p)-1).d _{p}}$, so we have 
$(r(p) - 1).d _{p}$ $\leq r(p)$, where $d _{p}$ is the degree of 
$\mathcal{G}(E (p)/E)$. In our case, this implies $d _{p} = 1$, 
which proves Corollary \ref{coro5.3}.
\end{proof}

\begin{rema}
\label{rema5.4}
Statement (5.1) and Corollary \ref{coro5.3} show that the equivalence 
(ii)$\leftrightarrow $(iv) in Proposition \ref{prop5.1} generalizes 
the concluding assertion of \cite{LLMS}, Theorem~1 (independently of 
the elementary type conjecture formulated in \cite{LLMS}).
\end{rema}

\medskip
Corollary \ref{coro5.3} gives us the possibility to determine the 
structure of the continuous character group $C(E(p)/E)$ of 
$\mathcal{G}(E(p)/E)$, for a $p$-quasilocal nonreal field $E$ 
containing a primitive $p$-th root of unity. Recall that $C(E(p)/E)$ 
is an abelian torsion $p$-group, whence, it decomposes into the 
direct sum $D(E(p)/E) \oplus R(E(p)/E)$, where $D(E(p)/E)$ is the 
maximal divisible subgroup of $C(E(p)/E)$ and $R(E(p)/E)$ is a 
(reduced) subgroup of $C(E(p)/E)$ isomorphic to 
$C(E(p)/E)/D(E(p)/E)$ (see \cite{Ka}, Ch. 7, Sect. 5, and \cite{F},  
Theorem~24.5). With this notation, our next result can be stated as 
follows:

\medskip
\begin{prop}
\label{prop5.5} 
Let $E$ be a $p$-quasilocal nonreal field, $\mu _{p}(E)$ the group 
of roots of unity in $E$ of $p$-primary degrees, and $r _{p}(E)$ the 
rank of $\mathcal{G}(E(p)/E)$ as a pro-$p$-group. Suppose that $\mu 
_{p}(E) \neq \{1\}$ and $\varepsilon _{p} \in E$ is a primitive 
$p$-th root of unity. Then:
\par
{\rm (a)} $C(E(p)/E) = D(E(p)/E)$ if and only if $\mu _{p}(E)$ is 
infinite or {\rm Br}$(E) _{p} = \{0\}$; when {\rm Br}$(E) _{p} \neq 
\{0\}$ and $\mu _{p}(E)$ is finite of order $p ^{\nu }$, the group 
$R(E(p)/E)$ is isomorphic to the maximal subgroup of {\rm Br}$(E)$ of 
period $p ^{\nu }$;  
\par
{\rm (b)} {\rm Br}$(E) _{p}$ is embeddable as a subgroup of 
$D(E(p)/E)$.
\end{prop}

\medskip
\begin{proof}
(a): It follows from Kummer theory that $C(E(p)/E) = D(E(p)/E)$, 
provided that $\mu _{p}(E)$ is infinite. We show that the same 
equality holds in the case of Br$(E) _{p} = \{0\}$. Then it follows 
from \cite{P}, Sect. 15.1, Proposition~b, that $\varepsilon _{p}$ 
lies in the norm group $N(L ^{\prime }/E)$, for every cyclic 
extension $L ^{\prime }$ of $E$ in $E(p)$; hence, by Albert's height 
theorem (cf. \cite{A1}, Ch. IX, Sect. 6, and \cite{FSS}, Sect. 2), 
there is a cyclic extension $L _{1} ^{\prime }$ of $E$ in $E(p)$, 
such that $L ^{\prime } \in I(L _{1} ^{\prime }/E)$ and $[L _{1} 
^{\prime }\colon L ^{\prime }] = p$. This implies $L ^{\prime } \in 
I(L _{1}/E)$, for some $\mathbb Z _{p}$-extension $L _{1}$ of $E$ in 
$E(p)$, and so proves that $C(E(p)/E) = D(E(p)/E)$. It remains to 
consider the case where Br$(E) _{p} \neq \{0\}$ and $\mu _{p}(E)$ 
has finite order $p ^{\nu }$. As $\mu _{p}(E) \neq \{1\}$, then we 
have $E(p) \neq E$, so it follows from \cite{Wh2}, Theorem~2, and the 
condition that $E$ is a nonreal field that $E(p)$ contains as a 
subfield a $\mathbb Z _{p}$-extension of $E$. In view of Galois 
theory, this means that $C(E(p)/E)$ possesses a quasicyclic 
$p$-subgroup, which proves that $D(E(p)/E) \neq \{0\}$. Note further 
that Br$(E) _{p}$ is a divisible group, since $E$ is nonreal and 
$p$-quasilocal (cf. \cite{Ch6}, I, Theorem~3.1). Therefore, by 
\cite{F}, Theorem~23.1, Br$(E) _{p}$ decomposes into the direct 
sum $\mathbb Z(p ^{\infty }) ^{d(p)}$ of isomorphic copies of the 
quasicyclic $p$-group $\mathbb Z(p ^{\infty })$, indexed by a set 
$I$ of cardinality $d(p)$ equal to the dimension of $_{p}{\rm 
Br}(E) \cong H ^{2}(\mathcal{G}(E(p)/E), \mathbb F _{p})$ as an 
$\mathbb F _{p}$-vector space. Moreover, it becomes 
clear that, for each $m \in \mathbb N$, the maximal subgroup of 
Br$(E) _{p}$ of period $p ^{m}$ decomposes into a direct sum of 
cyclic groups of order $p ^{m}$, indexed by $I$. Thus the latter 
conclusion of Proposition \ref{prop5.5} (a) is equivalent to the 
former part of the following assertions:
\par
\medskip
(5.4) (a) $C(E(p)/E)$ and the direct sum $D(E(p)/E) \oplus \mu 
_{p}(E) ^{d(p)}$ are isomorphic, where $\mu _{p}(E) ^{d(p)}$ is a 
direct sum of isomorphic copies of $\mu _{p}(E)$, indexed by a set 
of cardinality $d(p)$ (for a proof, see \cite{Ch1}, II, Lemma~2.3). 
\par
(b) A cyclic extension $M$ of $E$ in $E(p)$ is a subfield of a 
$\mathbb Z _{p}$-extension of $E$ in $E(p)$ if and only if there is 
$M ^{\prime } \in I(E(p)/M)$, such that $M ^{\prime }/E$ is cyclic 
and $[M ^{\prime }\colon M] = p ^{\nu }$; this is the case if and 
only if $\mu _{p}(E) \subset N(M/E)$.
\par
\medskip\noindent
The former part of (5.4) (b) is implied by (5.4) (a) and Galois 
theory, and the latter one follows from Albert's height theorem.
\par
(b): Proposition \ref{prop5.5} (a) allows us to consider only the 
case where Br$(E) _{p} \neq \{0\}$ and $\mu _{p}(E)$ has order $p 
^{\nu }$, for some $\nu \in \mathbb N$. Then it follows from (5.4)  
(a) and the nontriviality of $D(E(p)/E)$ that $r _{p}(E) \ge 2$. 
Using the observations preceding the statement of (5.4), one also 
sees that it is sufficient to prove the embeddability of $_{p}{\rm 
Br}(E)$ in $D(E(p)/E)$. Let now $\delta _{\nu }$ be a primitive $p 
^{\nu }$-th root of unity, and $M _{\lambda }$ be an extension of 
$E$ generated by a $p$-th root $\eta _{\lambda } \in E(p)$ of an 
element $\lambda \in E ^{\ast } \setminus E ^{\ast p}$. Then $M 
_{\lambda }/E$ is cyclic, $[M _{\lambda }\colon E] = p$ and 
$\mathcal{G}(M _{\lambda }/E)$ contains a generator $\sigma 
_{\lambda }$, such that the cyclic $E$-algebra $(M _{\lambda }/E, 
\sigma _{\lambda }, \delta _{\nu })$ is isomorphic to the symbol 
$E$-algebra $A _{\varepsilon _{p}}(\lambda , \delta _{\nu }; E)$. It 
is well-known that $A _{\varepsilon _{p}}(\lambda , \delta _{\nu }; 
E)$ and $A _{\varepsilon _{p}}(\delta _{\nu }, \lambda ; E)$ are 
inversely-isomorphic $E$-algebras. This, combined with \cite{P}, 
Sect. 15.1, Proposition~b, implies $\delta _{\nu } \in N(M _{\lambda 
}/E)$ if and only if $\lambda \in N(M _{\delta _{\nu }}/E)$. Hence, 
by (5.4) (b), the divisibility of Br$(E) _{p}$, and Theorem~23.1 of 
\cite{F}, Proposition \ref{prop5.5} (b) and the latter assertion of 
Proposition \ref{prop5.5} (a) are equivalent to the statement that 
$_{p}{\rm Br}(E)$ embeds as a subgroup of $N(M _{\delta _{\nu }}/E)/E 
^{\ast p}$. Since $N(M _{\delta _{\nu }}/E)/E ^{\ast p}$ is an 
abelian group of period $p$, this amounts to proving that $_{p}{\rm 
Br}(E)$ is its homomorphic image. We show that $_{p}{\rm Br}(E)$ is a 
homomorphic image of $N(M _{\mu }/E)/E ^{\ast p}$, for an arbitrary 
element $\mu \in E ^{\ast } \setminus E ^{\ast p}$. Fix $\mu 
^{\prime } \in E ^{\ast } \setminus E ^{\ast p}$ so that $M _{\mu '} 
\neq M _{\mu }$. Then $E ^{\ast }/N(M _{\mu '}/E) \cong $ $_{p}{\rm 
Br}(E)$, by \cite{P}, Sect. 15.1, Proposition~b (and the 
$p$-quasilocality of $E$), and $N(M _{\mu }/E)N(M _{\mu '}/E) = E 
^{\ast }$, by \cite{Ch6}, I, Lemma~4.3. Since, by Theorem~3.1, $N(M 
_{\mu }/E) \cap N(M _{\mu '}/E) = N(M _{\mu }M _{\mu '}/E)$, this 
yields
\par\vskip0.2truecm\noindent
$E ^{\ast }/N(M _{\mu '}/E) \cong N(M _{\mu }/E)/N(M _{\mu }M 
_{\mu '}/E)$, $E ^{\ast p} \le N(M _{\mu }M _{\mu '}/E)$ and
\par\noindent
$$N(M _{\mu }/E)/N(M 
_{\mu }M _{\mu '}/E) \cong (N(M _{\mu }/E)/E ^{\ast p})/(N(M _{\mu 
}M _{\mu '}/E)/E ^{\ast p});$$
\par\noindent
in particular, $_{p}{\rm Br}(E)$ is a homomorphic image of $N(M _{\mu 
}/E)/E ^{\ast p}$, which completes the proof of Proposition 
\ref{prop5.5} (b).
\end{proof}
\par
\medskip
The concluding results of this Section present applications of 
Proposition \ref{prop5.5} and Corollary \ref{coro5.3} to the study of 
$p$-primary index-exponent $K$-pairs, for a Henselian field $(K, v)$  
with a $p$-quasilocal field $\widehat K$, for some $p \in \mathbb P$. 
Recall first that, for any field $E$ and each $D \in d(E)$, exp$(D)$ 
divides ind$(D)$ and is divisible by any $p \in \mathbb P$ dividing 
ind$(D)$; in addition, $({\rm ind}(D), {\rm exp}(D))$is obtained as a 
componentwise product of index-exponent $l$-pairs, where $l$ runs 
across the set of prime divisors of ind$(D)$ (see \cite{P}, Sect. 
14.4). Index-exponent relations of algebras from $d(E)$ depend 
essentially on specific properties of $E$, and their description 
reduces to the special case of algebras $D _{p} \in d(E)$ of 
$p$-primary dimensions, for an arbitrary $p \in \mathbb P$. The 
study of index-exponent $p$-primary $E$-pairs relies on the 
knowledge of the Brauer $p$-dimension Brd$_{p}(E)$, defined as the 
least integer $b(p) \ge 0$, for which ind$(D _{p})$ divides exp$(D 
_{p}) ^{b(p)}$ whenever $D _{p} \in d(E)$ and $[D _{p}] \in {\rm 
Br}(E) _{p}$; when no such $b(p)$ exists, we say that Brd$_{p}(E)$ 
is infinite. Note that $(p ^{k}, p ^{n})\colon k, n \in \mathbb N, k 
\ge n$, are index-exponent $L$-pairs whenever $(L, \lambda )$ is a 
Henselian field, such that the quotient group $\lambda (L)/p\lambda 
(L)$ of the value group $\lambda (L)$ is infinite, $\widehat L$ is a 
nonreal field and $\mu _{p}(\widehat L) \neq \{1\}$. This is 
demonstrated by the proof of \cite{Ch8}, Corollary~4.5, which also 
shows that the result does not change, if the condition on $\mu 
_{p}(\widehat L)$ is replaced by the one that the rank $r 
_{p}(\widehat L)$ of $\mathcal{G}(\widehat L(p)/\widehat L)$ is 
infinite. Therefore, we restrict our considerations to the case 
where $\widehat K$ is a $p$-quasilocal nonreal field, $\mu 
_{p}(\widehat K) \neq \{1\}$ and $v(K)/pv(K)$ is a nontrivial finite 
group. Then Brd$_{p}(K)$ is determined by \cite{Ch12}, Theorem~4.1, 
as follows:
\par
\medskip
(5.5) If $v(K)/pv(K)$ has order $p ^{\tau (p)}$, then Brd$_{p}(K) = 
[(m _{p}(\widehat K) + \tau (p))/2]$, where $m _{p}(\widehat K) = 
{\rm min}\{r _{p}(\widehat K), \tau (p)\}$.
\par
\medskip\noindent
Relying on Proposition \ref{prop5.5}, we first show that if $r 
_{p}(\widehat K) = \infty $, then index-exponent $p$-primary 
$K$-pairs are fully determined by Brd$_{p}(K)$.
\par
\medskip
\begin{coro}
\label{coro5.6}
Let $(K, v)$ be a Henselian field with $\widehat K$ satisfying the 
conditions of Proposition \ref{prop5.5}, for some $p \in \mathbb P$. 
Assume also that $r _{p}(\widehat K) < \infty $, $v(K) \neq pv(K)$ 
and $v(K)/pv(K)$ is finite of order $p ^{\tau (p)}$. Then:
\par
{\rm (a)} $\widehat K$ has Galois extensions $U _{n}$, $U _{n} 
^{\prime }$ in $E(p)$, $n \in \mathbb N$, such that $[U 
_{n}\colon E] = p$, $[U _{1} \dots U _{n}\colon E] = p ^{n}$, 
$\mathcal{G}(U _{n} ^{\prime }/E) \cong \mathbb Z _{p}$ and $U _{n} 
\in I(U _{n} ^{\prime }/E)$, for each $n$;
\par
{\rm (b)} $(p ^{k}, p ^{n})\colon k, n \in \mathbb N, n \le k 
\le n\tau (p)$, are all nontrivial $p$-primary index-exponent pairs 
over $K$.
\end{coro}
\par
\medskip
\begin{proof}
(a): In view of (5.4) (b) and Galois theory, it suffices to prove 
that the maximal subgroup $_{p}D(\widehat K(p)/\widehat K)$ of 
$D(\widehat K(p)/\widehat K)$ of period $p$ is infinite. The infinity 
of $_{p}D(\widehat K(p)/\widehat K)$ follows from Proposition 
\ref{prop5.5} (b), if $_{p}{\rm Br}(\widehat K)$ is infinite. 
Note also that $r _{p}(\widehat K) = \infty $ if and only if 
$\widehat K ^{\ast }/\widehat K ^{\ast p}$ is infinite (cf. 
\cite{S2}, Ch. I, 4.1), which holds if and only if $C(\widehat 
K(p)/\widehat K)$ contains infinitely many elements of order $p$. It 
is therefore clear from (5.4) (a) that if $_{p}{\rm Br}(\widehat K)$ 
is finite, then $_{p}D(\widehat K(p)/\widehat K)$ is infinite, as 
required.
\par
(b): It follows from Corollary \ref{coro5.6} (a) and Galois theory
that, for each finite abelian $p$-group $G$, there exists a Galois
extension $U _{G}$ of $K$ in $K _{\rm ur}$ with $\mathcal{G}(U
_{G}/K) \cong G$. When the rank of $G$ is $\le \tau (p)$, one obtains 
from \cite{Mo}, Theorem~1 (or \cite{JW}, Example~4.3), that there is 
a nicely semi-ramified (abbr, NSR) algebra $N _{G} \in d(K)$, in the 
sense of \cite{JW}, possessing a maximal subfield $K$-isomorphic to 
$U _{G}$; this ensures that ind$(N _{G})$ and exp$(N _{G})$ equal the 
order and the period of $G$, respectively. Thus it becomes clear that 
there exist $N _{k,n} \in d(K)\colon k, n \in \mathbb N, n \le k \le 
\tau (p)n$, such that $N _{k,n}/K$ is NSR, ind$(N _{k,n}) = p ^{k}$ 
and exp$(N _{k,n}) = p ^{n}$. This proves Corollary \ref{coro5.6} 
(b), since by \cite{Ch12}, Theorem~4.1, the infinity of $r 
_{p}(\widehat K)$ implies $\tau (p) = {\rm Brd}_{p}(K)$.
\end{proof}

\medskip
Our next result supplements Corollary \ref{coro5.6} as follows:
\par
\medskip
\begin{coro}
\label{coro5.7}
Assume that $(K, v)$ is a Henselian field, such that $v(K)/pv(K)$ is 
finite of order $p ^{\tau (p)} > 1$, and $\widehat K$ is 
$p$-quasilocal with $\mu _{p}(\widehat K) \neq \{1\}$, $r 
_{p}(\widehat K) < \infty $, and $C(\widehat K(p)/\widehat K) = 
D(\widehat K(p)/\widehat K)$, for some $p \in \mathbb P$. Then  
nontrivial $p$-primary index-exponent $K$-pairs are described by as 
follows:
\par
{\rm (a)} $(p ^{k}, p ^{n})\colon k, n \in \mathbb N, n \le k \le 
n{\rm Brd}_{p}(K)$, if $\mu _{p}(\widehat K)$ is infinite;
\par
{\rm (b)} $(p ^{k}, p ^{n})\colon k, n \in \mathbb N$, $n \le k \le 
nm _{p}(\widehat K) + {\rm min}\{\nu , n\}([(\tau (p) - m 
_{p}(\widehat K))/2]$, in case $\mu _{p}(\widehat K)$ is finite of 
order $p ^{\nu }$, and $m _{p}(\widehat K)$ is defined as in (5.5).
\end{coro}

\medskip
\begin{proof}
The description of index-exponent $p$-primary $K$-pairs is obtained 
in the case of Corollary \ref{coro5.7} (b) by the method of proving 
\cite{Ch12}, Lemma~5.1 (replacing min$\{r _{p}(\widehat K) - 1, \tau 
(p)\}$ by $m _{p}(\widehat K)$). Henceforth, we assume that $\mu 
_{p}(\widehat K)$ is infinite. Note that, for any finite abelian 
$p$-group $G$ of rank $\le m _{p}(\widehat K)$, there exist a Galois 
extension $U _{G}$ of $K$ in $K _{\rm ur}$, and an algebra $N _{G} 
\in d(K)$, such that $\mathcal{G}(U _{G}/K) \cong G$, $N _{G}/K$ is 
NSR and $U _{G}$ is $K$-isomorphic to a maximal subfield of $N _{G}$. 
As in the proof of Corollary \ref{coro5.6} (b), one obtains further 
that there exist NSR-algebras $N _{k,n} \in d(K)\colon k, n \in 
\mathbb N$, $n \le k \le m _{p}(\widehat K)n$, with ind$(N _{k,n}) = 
p ^{k}$ and exp$(N _{k,n}) = p ^{n}$, for each admissible pair $k, 
n$. Since Brd$_{p}(K) = [(m _{p}(\widehat K) + \tau (p))/2]$, this 
fact proves Corollary \ref{coro5.7} in case $r _{p}(\widehat K) \ge 
\tau (p) - 1$. Suppose now that $\tau (p) - m _{p}(\widehat K) \ge 
2$, put $\bar m = {\rm Brd}_{p}(K) = [(m _{p}(\widehat K) + \tau 
(p))/2]$, fix a divisible hull $v(N _{G}) ^{\prime }$ of $v(N _{G})$, 
and take a finite abelian $p$-group $H$ of rank $r(H) \le [(\tau (p) 
- m _{p}(\widehat K))/2]$. Using \cite{Mo}, Theorem~1, and the 
natural bijection between $I(Y/K)$ and the set of subgroups of 
$v(Y)/v(K)$, for any totally ramified finite abelian $p$-extension 
$Y/K$ (cf. \cite{Sch}, Ch. 3, Theorem~2), one obtains that there is 
$T _{H} \in d(K)$ with the following properties: $T _{H}/K$ is totally 
ramified, $v(T _{H}) \subset v(N _{G}) ^{\prime }$ and $v(T 
_{H})/v(K)$ is isomorphic to $H \oplus H$; $v(N _{G}) \cap v(T _{H}) 
= v(K)$ and $T _{H} \otimes _{K} K _{\rm ur} \in d(K _{\rm ur})$; $T 
_{H}$ is a tensor product of $r(H)$ cyclic totally ramified 
$K$-algebras; ind$(T _{H})$ and exp$(T _{H})$ equal the order and the 
period of $H/v(K)$, respectively (see, e.g., \cite{JW}, for more 
details). Note also that, by \cite{Mo}, Theorem~1, $N _{G} \otimes 
_{K} T _{G} \in d(K)$. These observations indicate that, for any $n 
\in \mathbb N$, there exist $\Delta _{n} \in d(K)$ and $T _{n,\rho } 
\in d(K)\colon \rho = 1, \dots n \theta $, where $\theta = [(\tau (p) 
- m _{p}(\widehat K))/2]$, satisfying the following conditions:
\par
\medskip
(5.6) (a) $\Delta _{n}/K$ is NSR, exp$(\Delta _{n}) = p ^{n}$ and 
ind$(\Delta _{n}) = p ^{n.\bar m}$;
\par
(b) For each index $\rho $, $T _{n,\rho }/K$ is totally ramified, 
$\Delta _{n} \otimes _{K} T _{n,\rho } \in d(K)$, exp$(\Delta _{n} 
\otimes _{K} T _{n,\rho }) = p ^{n}$, and ind$(T _{n,\rho }) = p 
^{\rho }$.
\par
\medskip\noindent
Statements (5.6) (b) show that ind$(\Delta _{n} \otimes _{K} T 
_{n,\rho }) = p ^{n+\rho }$, $\rho = 1, \dots , \theta $, which 
completes the proof of Corollary \ref{coro5.7}. 
\end{proof}

\medskip
Summing-up Corollaries \ref{coro5.6}, \ref{coro5.7} and \cite{Ch12}, 
Lemma~5.1, and combining the latter with Corollary \ref{coro5.3}, one 
one fully describes $p$-primary index-exponent $K$-pairs, for a 
Henselian field $(K, v)$, such that $v(K) \neq pv(K)$ and $\widehat K$ 
is nonreal and $p$-quasilocal with $\mu _{p}(\widehat K) \neq \{1\}$. 
We refer the reader to \cite{Ch12}, Corollary~2.2 and Remark~4.2, for 
an analogous description in the case where $(K, v)$ is Henselian, $p 
= 2$, and $\widehat K$ is formally real $2$-quasilocal.

\medskip
\section{Preparation for the proof of Theorems 3.2 and 1.3}

\medskip
Let $E$ be a field, $M \in \Omega (E)$, $\Pi $ the set of prime
divisors of $[M\colon E]$, and $M _{p} = M \cap E (p)$, for each $p
\in \Pi $. When $M \neq E$, the homomorphism group ${\rm Hom}
(E ^{\ast }, \mathcal{G}(M/E))$ is isomorphic to the direct group
product $\prod _{p \in \Pi } {\rm Hom} (E ^{\ast },$ $\mathcal{G}(M
_{p}/E))$, so Theorem \ref{theo1.3} can be deduced from Theorem 
\ref{theo3.2}, Lemma \ref{lemm2.1} and the primary tensor 
decomposition of cyclic $E$-algebras (cf. \cite{P}, Sect. 15.3). The 
results of this Section serve as a basis for the proof of Theorem 
\ref{theo3.2}, presented in Section 7. Our starting point are the 
following two lemmas.

\medskip
\begin{lemm}
\label{lemm6.1} Let $E$ and $M$ be fields, such that $M \in \Omega
_{p} (E)$, for some $p \in P(E)$. Suppose that  $\mathcal{G}(M/E)$
has rank $t \geq 1$ as a $p$-group, and $F$ is an intermediate
field of $M/E$ of degree $[F\colon E] = p$. Then there exist cyclic
extensions $E _{1}, \dots , E _{t}$ of $E$ in $M$ with $E _{1} \dots
E _{t} = M$, $\prod _{i=1} ^{t} [E _{i}\colon E] = [M\colon E]$ and
$\prod _{i=1} ^{t} [(FE _{i})\colon F] = [M\colon F]$.
\end{lemm}

\begin{proof} In view of Galois theory, this is equivalent to the
following statement:

\medskip
(6.1) Let $G$ be a finite abelian $p$-group of rank $t \geq 1$ and
$H$ a maximal subgroup of $G$. Then these exist cyclic subgroups $G
_{1}, \dots , G _{t}$ of $G$, such that the (inner) products of
$G _{i}\colon \ i = 1,\dots , t$, and $H \cap G _{i}\colon \ i = 1,
\dots , t$, are direct and equal to $G$ and $H$, respectively.

\medskip
To prove (6.1) take a cyclic subgroup $G _{1} \le G$ of maximal
order so that the order of the group $H _{1} = H \cap G _{1}$ equals
the exponent of $H$. The choice of $G _{1}$ ensures that $G _{1}G
_{0} = G$ whenever $G _{0}$ is maximal among the subgroups of $G$,
which trivially intersect $G _{1}$ (cf. \cite{F}, Sects. 15 and 27).
This implies the existence of subgroups $H _{1} ^{\prime } \le H$
and $G _{1} ^{\prime } \le G$, such that $H _{1} \cap H _{1}
^{\prime } = G _{1} \cap G _{1} ^{\prime } = \{1\}$, $H _{1}
^{\prime } \le G _{1} ^{\prime }$ and the products $H _{1}H _{1}
^{\prime }$ and $G _{1}G _{1} ^{\prime }$ are equal to $H$ and $G$,
respectively. Therefore, (6.1) can be proved by induction on $t$.
\end{proof}

\medskip
\begin{lemm}
\label{lemm6.2}
Assume that $E$, $M$, $p$ and $t$ satisfy the conditions of Lemma
\ref{lemm6.1}, and let $F$ be a maximal subfield of $M$ including $E$. Then
there exist cyclic extensions $E _{1}, \dots , E _{t}$ of $E$ in
$M$, such that $\prod _{i=1} ^{t} [E _{i}\colon E] = [M\colon E]$,
$\prod _{i=1} ^{t} [(E _{i} \cap F)\colon E] = [F\colon E]$, and
the compositums $E _{1} \dots E _{t}$ and $(E _{i} \cap F) \dots (E
_{t} \cap F)$ are equal to $M$ and $F$, respectively.
\end{lemm}

\begin{proof}
This is equivalent to the following statement:

\medskip
(6.2) Let $G$ be a finite abelian $p$-group of rank $t \geq 1$, $H$ a
subgroup of $G$ of order $p$, and $\pi $ the natural homomorphism of
$G$ on $G/H$. Then $G$ contains elements $g _{1}, \dots , g _{t}$,
for which the products of the cyclic groups $\langle g _{1}\rangle
, \dots , \langle g _{t}\rangle $ and $\langle \pi (g _{1})\rangle
, \dots , \langle \pi (g _{t})\rangle $ are direct and equal to $G$
and $G/H$, respectively.

\medskip
For the proof of (6.2), consider a cyclic subgroup $C _{1} \le G$
of maximal order, and such that $C _{1} \cap H$ is of minimal
possible order. Then one can find $C _{1} ^{\prime } \le G$ so that
$C _{1} \cap C _{1} ^{\prime } = \{1\}$, $C _{1}C _{1} ^{\prime } =
G$ and $H \subset C _{1} \cup C _{1} ^{\prime }$. This implies $G/H$
is isomorphic to the direct product $C _{1}H/H \times C _{1}
^{\prime }H/H$, which allows one to prove (6.2) by induction on $t$.
\end{proof}

\medskip
Let now $G$ be a finite abelian $p$-group of rank $t \geq 2$. A
subset $g = \{g _{1}, \dots , g _{t}\}$ of $G$ is called an ordered
basis of $G$, if $g$ is a basis of $G$ (i.e. the group product
$\langle g _{1}\rangle \dots \langle g _{t}\rangle $ is direct and
equal to $G$) and the orders of the elements of $g$ satisfy the
inequalities $o(g _{1}) \geq \dots \geq o(g _{t})$. It is known that
the automorphism group Aut$(G)$ acts transitively on the set
Ob$(G)$ of ordered bases of $G$. Arguing by induction on $t$, and
using the fact that cyclic subgroups of $G$ of order $o(g _{1})$ are
direct summands in $G$ (see \cite{F},  Sect. 27), one obtains that
Aut$(G)$ has the following system of generators:

\medskip
(6.3) Aut$(G) = \langle d(k, m; h), t(i, j, s; h)\colon \ h \in {\rm
Ob}(G)\rangle $, where $i, j, k, s$ and $m$ are integers with $1 \leq
k, i, j \leq t$, $i \neq j$, $1 < m < o(h _{k})$, $1 \le s < o(h
_{i})$, $\gcd {(m, o(h _{k}))} = 1$ and max$\{1, o(h _{i})/o(h
_{j})\} \mid s$, and $d(k, m; h)$, $t(i, j, s; h)$ are defined by
the data $d(k, m; h) (h _{k}) = h _{k} ^{m}$, $d(k,m; h) (h _{k'}) =
h _{k'}\colon \ k ^{\prime } \neq k$, and $t(i, j, s; h) (h
_{j}) = h _{j}h _{i} ^{s}$, $t(i, j, s; h)$ $(h _{j'}) = h
_{j'}\colon \ j ^{\prime } \neq j$.

\medskip \noindent
This allows us to prove the following two lemmas without serious
technical difficulties.

\medskip
\begin{lemm}
\label{lemm6.3}
Let $E$, $M$, $p$ and $t$ satisfy the conditions of Lemma \ref{lemm6.1}, and
let $E$ be a $p$-quasilocal field. Suppose that $B$ a cyclic
subgroup of $ \ {\rm Br}(E) _{p}$ of order $o(B)$ divisible by
$[M\colon E]$, $b$ is a generator of $B$ and $B(M/E)$ is the group
of those $\b \in E ^{\ast }$, for which $[(L/E, \tau , \b )] \in B$
whenever $L/E$ is a cyclic extension, $L \subseteq M$ and $\langle
\tau \rangle =  \mathcal{G}(L/E)$. Assume also that $E _{1}, \dots ,
E _{t}$ are cyclic extensions of $E$ in $M$, such that $E _{1} \dots
E _{t} = M$ and $[M\colon E] = \prod _{j=1} ^{t} [E _{j}\colon E]$,
and for each index $j \leq t$, let $E _{j} ^{\prime }$ be the
compositum of the fields $E _{i}\colon \ i \neq j$, $\tau _{j}$ a
generator of $\mathcal{G}(E _{j}/E)$, and $\sigma _{j}$ the unique
$E _{j} ^{\prime }$-automorphism of $M$ extending $\tau _{j}$. Then
there exist a group homomorphism $\omega _{M/E,b}\colon \ B(M/E) \to
{\mathcal G}(M/E)$, and elements $c _{1}, \dots , c _{t}$ of $E
^{\ast }$ with the following properties:
\par
{\rm (i)} $c _{j} \in N(E _{j} ^{\prime }/E)$ and $[(E _{j}/E, \tau
_{j}, c _{j})] = [o(B)/[E _{j}\colon E]].b$, for each index $j$; the
co-set $c _{j}N(M/E)$ is uniquely determined by $b$ and $\tau _{j}$;
\par
{\rm (ii)} $\omega _{M/E,b}$ is the unique homomorphism of
$B(M/E)$ on $\mathcal{G}(M/E)$ mapping $c _{j}$ into $\sigma _{j}\colon
\ j = 1, \dots , t$, and with a kernel equal to $N(M/E)$;
\par
{\rm (iii)} $\omega _{M/E, b}$ does not depend on the choice of
the $t$-tuples $(E _{1}, \dots , E _{t})$ and $(\tau _{1}, \dots ,
\tau _{t})$; it induces a group isomorphism
$B(M/E)/N(M/E) \cong \mathcal{G}(M/E)$.
\end{lemm}

\medskip
\begin{proof}
If $t = 1$, our assertions can be deduced from Theorem \ref{theo3.1} 
and the general theory of cyclic algebras (see \cite{P}, Sect.
15.1). Henceforth, we assume that $t \geq 2$. Using consecutively
Galois theory and Theorem \ref{theo3.1}, one obtains that $N(E _{j}/E)N(E _{j}
^{\prime }/E) = E ^{\ast }\colon \ j = 1, \dots , t$, which implies
the existence of elements $c _{1}, \dots , c _{t}$ of $E ^{\ast }$
with the properties required by Lemma \ref{lemm6.3} (i). Denote by $T(M/E)$
the subgroup of $E ^{\ast }$ generated by the set $N(M/E) \cup \{c
_{1}, \dots , c _{t}\}$. It is easily verified that $\sigma _{1},
\dots , \sigma _{t}$ and $c _{1}N(M/E), \dots , c _{t}N(M/E)$ form
bases of the groups $\mathcal{G}(M/E)$ and $T(M/E)/$ $N(M/E)$,
respectively. Therefore, there exists a unique homomorphism $\omega
_{M/E, b}$ of $T(M/E)$ on $\mathcal{G}(M/E)$, such that $\omega _{M/E,
b} (c _{j}) = \sigma _{j}\colon \ j = 1, \dots , t$, and ${\rm
Ker}(\omega _{M/E, b}) = N(M/E)$. The mapping $\omega _{M/E,
b}$ is surjective, so it induces canonically an isomorphism of
$T(M/E)/N(M/E)$ on $\mathcal{G}(M/E)$. We aim at proving that $T(M/E) =
B(M/E)$. Assuming that $[M\colon E] = p ^{m}$, and proceeding by
induction on $m$, one obtains that this can be deduced from Lemma
\ref{lemm6.2} and \cite{P}, Sect. 15.1, Corollary~b, if $\omega _{M/E, b}$ has
the property claimed by the former part of Lemma \ref{lemm6.3} (iii). In
order to establish this property (the crucial point in our proof),
consider another basis $\sigma _{1} ^{\prime }, \dots , \sigma _{t}
^{\prime }$ of $\mathcal{G}(M/E)$, and for each $j \in \{1, \dots ,
t\}$, let $H _{j}$ be the subgroup of $\mathcal{G}(M/E)$ generated by
the elements $\sigma _{i} ^{\prime }\colon \ i \neq j$, $F _{j}$ the
fixed field of $H _{j}$, $\tau _{j} ^{\prime }$ the $E$-automorphism of
$F _{j}$ induced by $\sigma _{j} ^{\prime }$, and $F _{j} ^{\prime }$ the
compositum of all $F _{i}$ with $i \neq j$. It follows from Galois theory
that $F _{1} \dots F _{t} = M$, $\prod _{i=1} ^{t} [F _{i}\colon E] = [M\colon
E]$, $F _{j}/E$ is a cyclic extension, $\mathcal{G}(F _{j}/E) = \langle \tau
_{j} ^{\prime }\rangle $, and $\sigma _{j} ^{\prime }$ is the unique
$F _{j} ^{\prime }$-automorphism of $M$ extending $\tau _{j}
^{\prime }$ ($j = 1, \dots , t$). Since $\omega _{M/E, b}$ is
surjective, it maps some elements $c _{1} ^{\prime }, \dots , c _{t}
^{\prime }$ of $T(M/E)$ into $\sigma _{1} ^{\prime }, \dots , \s
_{t} ^{\prime }$, respectively. Applying Lemma \ref{lemm6.2} and \cite{P},
Sect. 15.1, Corollary~b, one concludes that the proof of Lemma \ref{lemm6.3}
will be complete, if we show that $c _{j} ^{\prime }\in N(F _{j}
^{\prime }/E)$ and $[(F _{j}/E, \tau _{j} ^{\prime }, c _{j}
^{\prime })] = [o(B)/[F _{j}\colon E]].b$, for every index $j$.
\par
It is clearly sufficient to consider the special case in which
$\{\sigma _{1}, \dots , \sigma _{t}\}$ and $\{\sigma _{1} ^{\prime
}, \dots , \sigma _{t} ^{\prime }\}$ are ordered bases of ${\mathcal
G}(M/E)$ (which implies $[E _{\rho }\colon E] = [F _{\rho }\colon
E]$, $\rho = 1, \dots , t$). In view of \cite{P}, Sect. 15.1,
Corollary~a, and (6.3), one may assume in addition that $\sigma _{j}
^{\prime } = \sigma _{j}\sigma _{i} ^{s}$, $c _{j} ^{\prime } = c
_{j}c _{i} ^{s}$, $\sigma _{u} ^{\prime } = \sigma _{u}$ and $c _{u}
^{\prime } = c _{u}\colon \ u \neq j$, for some pair $(i, j)$ of
different indices, and some integer $s$ satisfying the inequalities
$1 \le s < o(\sigma _{i})$ and divisible by max$ \{1, o(\sigma
_{i})/o(\s _{j})\}$. Then it follows from Galois theory that $F
_{i}F _{j} = E _{i}E _{j}$, $F _{w} = E _{w}, \tau _{w} ^{\prime } =
\tau _{w}\colon \ w \neq i$, and $F _{u} ^{\prime } = E _{u}
^{\prime }\colon \ u \neq j$. One also sees that if $t \geq 3$ and
$y \in (\{1, \dots ,t\} \setminus \{i, j\})$, then $c _{j} ^{\prime
} \in N(E _{y}/E)$. Since the cyclic $E$-algebras $(E _{j}/E, \tau
_{j}, c _{j})$ and $(E _{j}/E, \tau _{j}, c _{j} ^{\prime })$ are
isomorphic (see \cite{P}, Sect. 15.1, Proposition~b), and by Lemma
\ref{lemm4.1}, $N(F _{j} ^{\prime }/E) = \cap _{u \neq j} N(F _{u}/E)$, this
reduces the proof of Lemma \ref{lemm6.3} (ii) to the one of the assertions
that $c _{j} ^{\prime } \in N(F _{i}/E)$ and there is an
$E$-isomorphism $(E _{i}/E, \tau _{i}, c _{i}) \cong (F _{i}/E,
\tau _{i} ^{\prime }, c _{i})$. We first show that $(E _{i}/E, \tau
_{i}, c _{i}) \cong (F _{i}/E, \tau _{i} ^{\prime }, c _{i})$. The
assumptions on $E _{i}, F _{i}$ and $E _{j} = F _{j}$ guarantee that
the field $E _{i}E _{j} = F _{i}E _{j}$ is isomorphic as an
$E$-algebra to $E _{i} \otimes _{E} E _{j}$ and $F _{i} \otimes _{E}
E _{j}$ (cf. \cite{P}, Sect. 14.7, Lemma~b). It is therefore easy to
see from the general properties of tensor products (cf. \cite{P},
Sect. 9.2, Proposition~c) that $(E _{i}/E, \tau _{i}, c _{i})
\otimes _{E} (E _{j}/E, \tau _{j}, c _{j}) \cong (F _{i}/E, \tau
_{i} ^{\prime }, c _{i}) \otimes _{E} (E _{j}/E, \tau _{j}, c _{j}c
_{i} ^{s'})$ as $E$-algebras, where $s ^{\prime } = s.[o(\sigma
_{j})/o(\sigma _{i})]$. Since $s ^{\prime } \in \ZZ $ and $c _{i}
\in N(E _{j}/E)$, there exists an $E$-isomorphism $(E _{j}/E, \tau
_{j}, c _{j})$ $\cong (E _{j}/E, \tau _{j}, c _{j}c _{i} ^{s'})$.
The obtained results prove that $(E _{i}/E, \tau _{i}, c _{i})$ and
$(F _{i}/E,$ $\tau _{i} ^{\prime }, c _{i})$ are similar over $E$.
As $[E _{i}\colon E] = [F _{i}\colon E]$, we also have $[(E _{i}/E,
\tau _{i}, c _{i})\colon E] = [(F _{i}/E, \tau _{i} ^{\prime }, c
_{i})\colon E]$ $= [E _{i}\colon E] ^{2}$, so it turns out that $(E
_{i}/E, \tau _{i}, c _{i}) \cong (F _{i}/E,$ $\tau _{i} ^{\prime },
c _{i})$, as claimed.
\par
It remains to be seen that $c _{j} ^{\prime } \in N(F _{i}/E)$.
Suppose first that $F _{i} \cap E _{i} = E$, take elements $\b
_{i} \in E _{j}$, $\b _{j} \in E _{i}$, $\d \in {\rm Br}(E)$
and an algebra $D \in d(E)$ so that $N _{E} ^{E _{j}} (\b _{i}) =
c _{i}$, $N _{E} ^{E _{i}} (\b _{j}) = c _{j}$, $([E _{i}\colon
E].[E _{j}\colon E])\d = b$, and $o(B)\d = [D]$, and denote
by $\varphi _{i}$ and $\varphi _{j}$ the automorphisms of $E _{i}E
_{j}$ induced by $\sigma _{i}$ and $\sigma _{j}$, respectively. It
is easily verified that $[E _{j}\colon E].[D] = [(E _{i}/E, \tau
_{i}, c _{i})]$ and $\mathcal{G}(E _{i}E _{j}/E _{j}) = \langle \varphi
_{i}\rangle $. This, combined with the fact that $E _{i}E _{j} = F
_{i}E _{j}$, $E _{i} \cap E _{j} = E$ and $\varphi _{i}$ extends
$\tau _{i}$, enables one to deduce from (4.1) (iii), Proposition \ref{prop2.2}
and the RC-formula that $D \otimes _{E} E _{j}$ is similar to the
cyclic $E _{j}$-algebra $(E _{i}E _{j}/E _{j}, \varphi _{i}, \b
_{i})$. Since $[E _{i}\colon E].[D] = [(E _{j}/E, \tau _{j}, c
_{j})]$, $\mathcal{G}(E _{i}E _{j}/E _{i}) = \langle \varphi
_{j}\rangle $, and $\varphi _{j}$ is a prolongation of $\tau _{j}$,
it is analogously proved that $[D \otimes _{E} E _{i}] = [(E _{i}E
_{j}/E _{i}, \varphi _{j}, \b _{j})]$ in Br$(E _{i})$. At the
same time, it is clear from Galois theory and the condition $F _{i}
\cap E _{i} = E$ that $\mathcal{G}(E _{i}E _{j}/E) = \langle \varphi
_{j}, \varphi _{i} ^{s}\rangle $ (i.e. ${\rm g.c.d.}(s, p) = 1$),
$[E _{i}\colon E] \leq [E _{j}\colon E]$, and $E _{i}F _{i} = E
_{i}\Phi $, where $\Phi $ is the extension of $E$ in $E _{j}$ of
degree $[\Phi \colon E] = [E _{i}\colon E]$. Note also that
$\varphi _{j}$ extends $\tau _{i} ^{\prime -s}$, since $\sigma _{j}
^{\prime } (\lambda ) = \lambda $, for every $\lambda \in F _{i}$.
Observing now that $[(E _{i}/E, \tau _{i}, c _{i})] = \omega .[(E
_{j}/E, \tau _{j}, c _{j})]$ and $[(E _{i}F _{i}/E _{i}, \psi _{j},
\b _{j})] = \omega .[(E _{i}E _{j}/E _{i}, \varphi _{j}, \b
_{j})]$, where $\omega = [E _{j}\colon E]/[E _{i}\colon E] = [E
_{i}E _{j}\colon E _{i}F _{i}]$ and $\psi _{j}$ is the automorphism
of $E _{i}F _{i}$ induced by $\varphi _{j}$ (cf. \cite{P}, Sect.
15.1, Corollary~b), one obtains from (4.1) (iii) that $(E _{i}/E,
\tau _{i}, c _{i}) \cong (\Phi /E, \bar \psi _{j}, c _{j}) \cong (F
_{i}/E, \tau _{i} ^{\prime -s}, c _{j})$ over $E$ ($\bar \psi _{j}$
being the automorphism of $\Phi $ induced by $\varphi _{j}$). Since
$(E _{i}/E, \tau _{i}, c _{i}) \cong (F _{i}/E, \tau _{i} ^{\prime
}, c _{i})$, these results indicate that $(F _{i}/E, \tau _{i}
^{\prime -s}, c _{j}) \cong (F _{i}/E, \tau _{i} ^{\prime -s}, c
_{i} ^{-s})$, which means that $c _{j} ^{\prime } \in N(F _{i}/E)$.
\par
Let now $E _{i} \cap F _{i} = \widetilde E$, $[\widetilde E\colon E]
= \mu > 1$, $s/\mu = \tilde s$ and $E _{i}E _{j} =  E _{i,j}$. Then
$E _{i,j} \in \Omega (\widetilde E)$ and $\mathcal{G}(E
_{i,j}/\widetilde E) = \mathcal{G}(E _{i,j}/E _{i})\mathcal{G}(E
_{i,j}/F _{i})$, which implies $\tilde s \in \ZZ $ and ${\rm
g.c.d.}(\tilde s, p) = 1$. For each index $u \neq i$, take $\g
_{u} \in E _{u} ^{\prime }$ so that $N _{E} ^{E _{u}'} (\g _{u})
= c _{u}$, and put $\tilde c _{u} = N _{\widetilde E} ^{E _{u}'}
(\g _{u})$, $\widetilde E _{u} = E _{u}\widetilde E$ and
$\widetilde F _{u} = F _{u}\widetilde E$. Since $E _{i} \cap E _{u}
= E$, it follows from (4.1) (ii) that $E _{i} \cap \tilde E
_{u} = \widetilde E$, $[\widetilde E _{u}\colon \widetilde E] = [E
_{u}\colon E]$, $\tau _{u}$ uniquely extends to a $\widetilde
E$-automorphism $\tilde \tau _{u}$ of $\widetilde E _{u}$, and
$\mathcal{G}(\widetilde E _{u}/\widetilde E) = \langle \tilde \tau
_{u}\rangle $. Note also that $N _{E} ^{E _{u}} (\theta _{u}) = N
_{\widetilde E} ^{\widetilde E _{u}} (\theta _{u})$, $\theta _{u}
\in E _{u}$, whence $c _{i} \in N(\widetilde E _{u}/\widetilde E)$.
Fix elements $b _{\mu } \in {\rm Br}(E)$ and $\Delta _{u} \in d(E)$
so that $\mu .b _{\mu } = b$ and $[\Delta _{u}] = [o(B)/[E _{u}\colon
E]]b _{\mu }$, and put $b ^{\prime } =  \rho _{E/\widetilde E} (b
_{\mu })$. It is clear from the double centralizer theorem (cf.
\cite{P}, Sects. 12.7 and 13.3) that $[(E _{i}/E, \tau _{i}, c _{i})
\otimes _{E} \widetilde E] = [(E _{i}/\widetilde E, \tau _{i} ^{\mu
}, c _{i})]$ in Br$(\widetilde E)$. Applying Proposition \ref{prop2.2}, (4.1)
(iii) and the RC-formula, one obtains consecutively that $o(B)$
equals the order of $b ^{\prime }$ in Br$(\widetilde E)$, $[(E
_{i}/E, \tau _{i}, c _{i})] = [o(B)/[E _{i}\colon \widetilde E]]b
_{\mu }$, $[(E _{i}/\widetilde E, \tau _{i} ^{\mu }, c _{i})] =
[o(B)/[E _{i}\colon \widetilde E]].b ^{\prime }$, and for each $u
\neq i$, $[(E _{u}/E, \tau _{u}, c _{u})] = \mu [\Delta _{u}]$ and
$[(\widetilde E _{u}/\widetilde E, \tilde \tau _{u}, \tilde c _{u})] =
[\Delta _{u} \otimes _{E} \widetilde E] = [o(B)/[\widetilde E
_{u}\colon \widetilde E]].b ^{\prime }$. Consider now
$M/\widetilde E$, $\tilde s$, $b ^{\prime }$, the fields $E _{i}, F
_{i}$, $\widetilde E _{u}, \widetilde F _{u}\colon \ u \neq i$, and
the algebras $(E _{i}/\widetilde E, \tau _{i} ^{\mu }, c _{i})$, $(F
_{i}/\widetilde E, \tau _{i} ^{\prime \mu }, c_{i})$, $(\widetilde E
_{j}/\widetilde E, \tilde \tau _{j}, \tilde c _{j})$, instead of $M/E$,
$s$, $b$, $E _{u'}, F _{u'}\colon \ u ^{\prime } = 1, \dots , t$, and
$(E _{i}/E, \tau _{i}, c _{i})$, $(F _{i}/E, \tau _{i} ^{\prime }, c
_{i})$, $(E _{j}/E, \tau _{j},$ $c _{j})$, respectively. As in the
proof of our assertion in case ${\rm g.c.d.}(s, p) = 1$, one
concludes that $\tilde c _{j}c _{i} ^{\tilde s} \in N(F _{i}/\widetilde
E)$. Since $N _{E} ^{\widetilde E} (\tilde c _{j}c _{i} ^{\tilde s}) = c
_{j} ^{\prime }$, this yields $c _{j} ^{\prime } \in N(F _{i}/E)$ (and
$T(M/E) = B(M/E)$), which completes the proof of Lemma \ref{lemm6.3}.
\end{proof}

\medskip
\begin{lemm}
\label{lemm6.4}
In the setting of Lemma \ref{lemm6.3}, let $F$ be an extension of $E$ in $M$,
$\pi _{M/F}$ the natural projection of $\mathcal{G}(M/E)$ on
${\mathcal G}(F/E)$, $\rho _{E/F} (b) = b ^{\prime }$, $B ^{\prime }
= \langle b ^{\prime }\rangle $, and $B _{\kappa } = \langle b
_{\kappa }\rangle $, where $b _{\kappa } = \kappa .b$, for some
$\kappa \in \mathbb N$ dividing $o(B)/[M\colon E]$. Then $[M\colon
E] \mid o(B _{\kappa })$, $[M\colon F] \mid o(B ^{\prime })$, and
the mappings $\omega _{M/E, b}$, $\omega _{M/E, b _{\kappa }}$,
$\omega _{F/E, b}$ and $\omega _{M/F, b'}$, determined as required
by Lemma \ref{lemm6.3}, are related as follows:
\par
{\rm (i)} $\omega _{M/E, b} = \omega _{M/E, b _{\kappa }}$ and
$\omega _{F/E, b} =  \pi _{M/F} \circ \omega _{M/E, b}$;
\par
{\rm (ii)} $\omega _{M/F, b'} (\g ) = \omega _{M/E, b} (N _{E}
^{F} (\g ))$, for every $\g \in F ^{\ast }$.
\end{lemm}

\begin{proof}
Clearly, $o(B _{\kappa }) = o(B)/\kappa $, and Proposition \ref{prop2.2}
implies that $o(B ^{\prime }) = o(B)/[F\colon E]$, so the assertions
that $[M\colon E] \mid o(B _{\kappa })$, $[M\colon F] \mid o(B
^{\prime })$ and $\omega _{M/E, b} = \omega _{M/E, b _{\kappa }}$
are obvious. Note also that if $t = 1$, then the remaining
statements of the lemma can be deduced from the general properties
of cyclic algebras. Suppose further that $t \geq 2$. First we prove
that $\omega _{F/E, b} = \pi _{M/F} \circ \omega _{M/E, b}$ in case
$F$ is a maximal subfield of $M$. Let $E _{1}, \dots , E _{t}$ be
extensions of $E$ in $M$ with the properties described by Lemma \ref{lemm6.2}.
Then there exists an index $j$, such that $E _{j} \cap F$ is a
maximal subfield of $E _{j}$ and $E _{i} \subseteq F$, for any other
index $i$. This allows us to prove our assertion by applying Lemma
\ref{lemm6.3} and \cite{P}, Sect. 15.1, Corollary~b.
\par
We turn to the proof of Lemma \ref{lemm6.4} (ii) in the special case of
$[F\colon E] = p$. It is briefly presented, since we argue along the
same lines as the concluding part of the proof of Lemma \ref{lemm6.3}. Take $E
_{1}, E _{1} ^{\prime }, \dots , E _{t}, E _{t} ^{\prime }$, $\tau
_{1}, \sigma _{1}, \dots , \tau _{t}, \sigma _{t}$ and $c _{1},
\dots , c _{t}$ in accordance with Lemmas \ref{lemm6.1} and \ref{lemm6.3}. Then $F
\subseteq E _{j}$ and $E _{u} \cap F = E\colon \ u \neq j$, for some
index $j$. This implies $F \subseteq E _{u} ^{\prime }$ and $\tau
_{u}$ uniquely extends to an $F$-automorphism $\tilde \tau _{u}$ of
the field $E _{u}F := \tilde E _{u}$, when $u$ runs across the set
$W _{j} = \{1, \dots , t\} \setminus \{j\}$. Choose elements $\alpha
_{1} \in E _{1} ^{\prime }, \dots , \alpha _{t} \in E _{t} ^{\prime
}$ so that $N _{E} ^{E _{k}'} (\alpha _{k}) = c _{k}\colon \ k = 1,
\dots ,t$, and put $\tilde E _{j} = E _{j}$, $\tilde c _{j} = c
_{j}$, $\tilde \tau _{j} = \tau _{j} ^{p}$, $\tilde \sigma _{j} =
\sigma _{j} ^{p}$, and $\tilde c _{u} = N _{F} ^{E _{u}'} (\alpha
_{u}), \tilde \sigma _{u} = \sigma _{u}\colon \ u \in W _{j}$. It is
easily verified that $N _{E} ^{F} (\tilde c _{j}) = c _{j} ^{p}$ and
$N _{E} ^{F} (\tilde c _{u}) = c _{u}$, $u \in W _{j}$. Therefore,
it follows from (4.1) (iii), \cite{P}, Sect. 14.7, Lemma~a, and the
equality $o(B ^{\prime }) = o(B)/p$ that $M/F$, $b ^{\prime }$, and
$\tilde E _{\rho }$, $\tilde \tau _{\rho }, \tilde \sigma _{\rho },
\tilde c _{\rho }$, where $\rho $ runs through $W _{j}$ or $\{1,
\dots , t\}$ depending on whether or not $F = E _{j}$, are related
as in Lemma \ref{lemm6.3}. Since the RC-formula and \cite{P}, Sect. 15.1,
Proposition~b, yield Cor$_{F/E} ([(E _{j}/F, \tau _{j} ^{p}, c
_{j})]) = [(E _{j}/E, \tau _{j}, c _{j} ^{p})]$, this proves Lemma
\ref{lemm6.4} (ii) when $[F\colon E] = p$.
\par
Let finally $F$ be any proper extension of $E$ in $M$ different
from $M$. By Galois theory and the structure of finite abelian
groups, then $M$ possesses subfields $F _{0}$ and $M _{0}$, such
that $E \subset F _{0} \subseteq F \subseteq M _{0}$ and $[F
_{0}\colon E] = [M\colon M _{0}] = p$. In view of the transitivity
of norm mappings, canonical projections of Galois groups and
scalar extension maps of Brauer groups in towers of intermediate
fields of $M/E$ (cf. \cite{P}, Sect. 9.4, Corollary~a), the
considered special cases of Lemma \ref{lemm6.4} enable one to complete
inductively its proof.
\end{proof}
\medskip

\section{Existence and form of Hasse symbols}

\medskip
The main purpose of this Section is to prove Theorem \ref{theo3.2}. Fix an
$\mathbb F _{p}$-basis $I _{p}$ of $_{p}{\rm Br}(E)$ and a
generator $\varphi _{\infty }$ of $\mathcal{G}(E _{\infty }/E)$ as a
topological group. Take a subset $\Lambda _{p} (E) = \{b _{i,n}
(p)\colon \ i \in I _{p}, n \in \mathbb N\}$ of Br$(E) _{p}$ so that
$\{b _{i,1} (p)\colon \ i \in I _{p}\}$ is a basis of $_{p} {\rm
Br}(E)$ as an $\mathbb F _{p}$-vector space, and $pb _{i,n} (p) = b
_{i,n-1} (p)\colon \ i \in I _{p}, n \in \ZZ $ and $n \geq 2$ (the
existence of $\Lambda _{p} (E)$ is implied by Proposition \ref{prop2.2} (ii)).
For each $E ^{\prime } \in \Omega _{p} (E)$, let $\rho _{E/E'} (b
_{i,n} (p)) = b _{i,n} (E ^{\prime })\colon \ (i, n) \in I _{p}
\times \mathbb N$, $E ^{\prime } _{\infty } = E _{\infty }E ^{\prime
}$, and $[(E ^{\prime } \cap E _{\infty })\colon E] = m(E ^{\prime
})$. As $\mathcal{G}(E _{\infty }/E) \cong \mathbb Z _{p} \cong H$
whenever $H$ is an open subgroup of $\mathbb Z _{p}$ (see \cite{S2},
Ch. I, 1.5, 4.2), it follows from Galois theory that $E _{\infty }
^{\prime }/E ^{\prime }$ is a $\mathbb Z _{p}$-extension and
$\varphi _{\infty } ^{m(E')}$ uniquely extends to an $E ^{\prime
}$-automorphism $\varphi _{\infty } (E ^{\prime })$ of $E _{\infty}
^{\prime}$; one also sees that $\varphi _{\infty } (E ^{\prime })$
topologically generates $\mathcal{G}(E _{\infty } ^{\prime }/E
^{\prime })$. For each $n \in \mathbb N$ and $i \in I _{p}$, denote
by $\varphi _{n} (E ^{\prime })$ the automorphism of $\Gamma _{n}E
^{\prime }$ induced by $\varphi _{\infty } (E ^{\prime })$, and by
$g _{i,n} (E ^{\prime })$ the element of $\mathcal{G}(\Gamma _{n}E
^{\prime }/E ^{\prime }) ^{d(p)}$ with components $g _{i,n} (E
^{\prime }) _{i} = \varphi _{n} (E ^{\prime })$, and $g _{i,n} (E
^{\prime }) _{i'} = 1\colon \ i ^{\prime } \in I _{p} \setminus
\{i\}$. By Proposition \ref{prop2.2} and \cite{P}, Sect. 15.1, Proposition~a,
$E ^{\prime \ast }$ has a subset $C _{p} (E ^{\prime }) = \{c _{i,n}
(E ^{\prime })\colon \ i \in I _{p}, n \in \mathbb N\}$, such that
$[((\Gamma _{n}E ^{\prime })/E ^{\prime }, \varphi _{n} (E ^{\prime
}), c _{i,n} (E ^{\prime }))] = b _{i,n} (E ^{\prime })$, for each
$(i, n) \in I _{p} \times \mathbb N$. Observe that, for any $n \in
\mathbb N$, there is a unique surjective homomorphism $( \ , (\Gamma
_{n}E ^{\prime })/E ^{\prime })\colon \ E ^{\prime \ast } \to
\mathcal{G}(\Gamma _{n}E ^{\prime }/E ^{\prime }) ^{d(p)}$, whose
kernel is $N(\Gamma _{n}E ^{\prime }/E ^{\prime })$, and which maps
$c _{i,n} (E ^{\prime })$ into $g _{i,n} (E ^{\prime })$, when $i
\in I _{p}$. Fix some $M ^{\prime } \in \Omega _{p} (E ^{\prime })$
and $\mu  \in \ZZ $ so that $p ^{\mu } \geq [M ^{\prime }\colon E]$,
and for each $i \in I _{p}$, let $B(M ^{\prime }/E ^{\prime }) _{i}
= \{\b _{i} \in E ^{\prime \ast }\colon \ [(L ^{\prime }/E ^{\prime
}, \sigma ^{\prime }, \b _{i})] \in \langle b _{i, \mu } (E ^{\prime
})\rangle \}$, for every cyclic extension $L ^{\prime }$ of $E
^{\prime }$ in $M ^{\prime }$, where $\sigma ^{\prime }$ is a
generator of $\mathcal{G}(L ^{\prime }/E ^{\prime })$. It follows
from Proposition \ref{prop2.2} (ii)-(iii), the definition of $\Lambda _{p}
(E)$ and Lemma \ref{lemm6.3} that the groups $\overline B(M ^{\prime }/E
^{\prime }) _{i} := B(M ^{\prime }/E ^{\prime }) _{i}/N(M ^{\prime
}/E ^{\prime })$, $i \in I _{p}$, have the following property:

\medskip
(7.1) $\overline B(M ^{\prime }/E ^{\prime }) \cong \mathcal{G}(M
^{\prime }/E ^{\prime })$, for each index $i$, and the inner product
of $\overline B(M ^{\prime }/E ^{\prime }) _{i}\colon \ i \in I
_{p}$, is direct and coincides with $E ^{\prime \ast }/N(M ^{\prime
}/E ^{\prime })$.

\par
\medskip \noindent
Therefore, there exists a unique homomorphism $( \ , M ^{\prime }/E
^{\prime })$ of $E ^{\prime \ast }$ into

\noindent
$\mathcal{G}(M ^{\prime }/E ^{\prime
}) ^{d(p)}$, mapping $B(M ^{\prime }/E ^{\prime }) _{i}$ into the
$i$-th component of $\mathcal{G}(M ^{\prime }/E ^{\prime }) ^{d(p)}$
by the formula $\b _{i} \to  \omega _{M'/E',b _{i,\mu } (E')} (\b
_{i})$, for each $i \in I _{p}$, where $\omega _{M'/E',b _{i,\mu }
(E')}$ is defined as in Lemma \ref{lemm6.3}. Hence, ${\rm Ker}( \ ,M ^{\prime
}/E ^{\prime }) = N(M ^{\prime }/E ^{\prime })$, and by Lemma \ref{lemm6.4},
the sets $H(E ^{\prime }) = \{( \ , M ^{\prime }/E ^{\prime })\colon
M ^{\prime } \in \Omega _{p} (E ^{\prime })\}$, $E ^{\prime } \in
\Omega _{p} (E)$, consist of surjections related as required by
Theorem \ref{theo3.2} (ii)$\div $(iii).

\par
Suppose now that $\Theta _{p} (E ^{\prime }) = \{\theta (M ^{\prime
}/E ^{\prime })\colon \ E ^{\prime \ast } \to \mathcal{G}(M ^{\prime
}/E ^{\prime }) ^{d(p)}, M ^{\prime } \in \Omega _{p} (E ^{\prime
})\}$, $E ^{\prime } \in \Omega _{p} (E)$, is a system of surjective
homomorphisms with the same kernels and relations, and such that
$\theta (\Gamma _{n}/E) = ( \ , \Gamma _{n}/E)$, for every $n \in
\mathbb N$. Then it follows from Proposition \ref{prop2.2}, the RC-formula and
(4.1) (iii) that $\theta (\Gamma _{n}E ^{\prime }/E ^{\prime }) = (
\ , (\Gamma _{n}E ^{\prime })/E ^{\prime })$, for each pair $(E
^{\prime }, n) \in \Omega _{p} (E) \times \mathbb N$. We show that
$\Theta _{p} (E ^{\prime }) = H _{p} (E ^{\prime })$, for any $E
^{\prime } \in \Omega _{p} (E)$. This is obvious, if $E _{\infty } =
E (p)$, so we assume further that $E (p) \neq E _{\infty }$. In view
of Proposition \ref{prop2.2} (i) and the already established part of Theorem
\ref{theo3.2}, it is sufficient to prove that $\theta (M/E) = ( \ , M/E)$, for
an arbitrary fixed field $M \in \Omega _{p} (E)$. Note first that
the compositum $ME _{\infty }$ possesses a subfield $M _{0} \in
\Omega _{p} (E)$, such that $M _{0} \cap E _{\infty } = E$ and $M
_{0}E _{\infty } = ME _{\infty }$. This follows from Galois theory
and the projectivity of $\mathbb Z _{p}$ as a profinite group (cf.
\cite{S2}, Ch. I, 5.9). In particular, $M \subseteq M _{0}\Gamma
_{n}$, for every sufficiently large index $n$. We also have $\theta
(M/E) = \pi _{M _{0}\Gamma _{n}/M} \circ \theta (M _{0}\Gamma
_{n}/E)$, which allows us to consider only the special case of $M =
M _{0}\Gamma _{\kappa }$ and $M _{0} \neq E$, where $\kappa $ is
chosen so that $[M _{0}\colon E] \mid p ^{\kappa }$. Let now $t$ be
the rank of $\mathcal{G}(M/E)$, and $E _{1}, \dots , E _{t}$ be
cyclic extensions of $E$ in $M$, such that $\prod _{u=1} ^{t} [E
_{u}\colon E] = [M\colon E]$, $E _{1} = \Gamma _{\kappa }$ and $E
_{2} \dots E _{t} = M _{0}$. Take $E _{1} ^{\prime }, \dots , E _{t}
^{\prime }$ as in Lemma \ref{lemm6.3}, denote for convenience by $\tau _{1}$
the automorphism $\varphi _{\kappa } (E)$ of $E _{1}$, and let $\tau
_{u}$ be a generator of $\mathcal{G}(E _{u}/E)$, for every $u \in
\{2, \dots , t\}$. Fix an index $x \in I _{p}$, put $b _{x, \kappa }
(p) = b _{x}$, and identifying $\mathcal{G}(M/E)$ with the $x$-th
component of $\mathcal{G}(M/E) ^{d(p)}$, consider elements $x _{1},
\dots , x _{t}$ of $E ^{\ast }$ determined so that $\theta (M/E) (x
_{u})$ equals the $E _{u} ^{\prime }$-automorphism $\sigma _{u}$ of
$M$ extending $\tau _{u}$, for each positive integer $u \leq t$. The
assumptions on $\Theta _{p} (E ^{\prime })\colon \ E ^{\prime } \in
\Omega _{p} (E)$, imply that $x _{u} \in N(E _{u} ^{\prime
}/E)\colon \ u = 1, \dots , t$, and $[(E _{1}/E, \tau _{1}, x _{1})]
= b _{x}$. We show as in the proof of Lemma \ref{lemm6.3} (iii) that $[(E
_{u}/E, \tau _{u}, x _{u})] = (p ^{\kappa }/[E _{u}\colon E])b
_{x}$, for each $u$. Fix an index $u \geq 2$, put $\omega _{u} = p
^{\kappa }/[E _{u}\colon E]$, $\sigma _{1} ^{\prime } = \sigma
_{1}\sigma _{u}$, $x _{1} ^{\prime } = x _{1}x _{u}$, and denote by
$F _{u}$ the fixed field of the subgroup of $\mathcal{G}(M/E)$
generated by $\sigma _{1} ^{\prime }$ and $\sigma _{u'}\colon \ u
^{\prime } \notin \{1, u\}$. It is easily verified that $F _{1}
\dots F _{t} = M$ and $\prod _{u'=1} ^{t} [F _{u'}\colon E] =
[M\colon E]$, where $F _{u'} = E _{u'}\colon \ u ^{\prime } \neq u$.
As ${\rm Ker} (\theta (F _{u}/E)) = N(F _{u}/E)$, $\theta (M/E) (x
_{1} ^{\prime }) = \sigma _{1} ^{\prime }$ and $\sigma _{1} ^{\prime
} \in \mathcal{G}(M/F _{u})$, the equality $\theta (F/E) = \pi
_{M/E} \circ \theta (M/E)$ ensures that $x _{1} ^{\prime } \in N(F
_{u}/E)$. Observe also that $E _{1}E _{u} = E _{1}F _{u}$ and $E
_{u} \cap F _{u} = E$. This implies that $(E _{1}/E, \tau _{1}, x
_{1}) \otimes _{E} (E _{u}/E, \tau _{u}, x _{u})$ is $E$-isomorphic
to $(E _{1}/E, \tau _{1} ^{\prime }, x _{1}x _{u} ^{\omega _{u}})
\otimes _{E} (F _{u}/E, f _{u}, x _{u})$, where $\tau _{1} ^{\prime
} \in \mathcal{G}(E _{1}/E)$ and $f _{u} \in {\mathcal G}(F _{u}/E)$
are induced by $\sigma _{1} ^{\prime }$ and $\sigma _{u}$,
respectively. Since $\sigma _{u} \in \mathcal{G}(M/E _{1})$ and $x
_{u} \in N(E _{1}/E)$, it thereby turns out that $(E _{u}/E, \tau
_{u}, x _{u}) \cong (F _{u}/E, f _{u}, x _{u}) \cong (F _{u}/E, f
_{u} ^{-1}, x _{1}) \cong (F _{u}/E, f _{1}, x _{1})$ over $E$, $f
_{1}$ being the automorphism of $F _{u}$ induced by $\sigma _{1}$.
Furthermore, if $\Phi _{u}$ is the extension of $E$ in $E _{1}$ of
degree $[E _{u}\colon E]$, then $E _{u}\Phi _{u} = E _{u}F _{u}$
and, because $x _{1} \in N(E _{u}/E)$, it follows that $[(E _{u}/E,
\tau _{u}, x _{u})] = [(F _{u}/E, f _{1}, x _{1})] = \omega _{u}[(E
_{1}/E, \tau _{1}, x _{1})] = \omega _{u}b _{x}$, as claimed. The
obtained result indicates that $\theta (M/E) (\alpha _{x}) = (\alpha
_{x}, M/E)$, for each $\a _{x} \in B(M/E) _{x}$. As $x$ is an
arbitrary element of $I _{p}$, this enables one to complete the
proof of Theorem \ref{theo3.2} by applying (7.1).

\medskip
\begin{rema}
\label{rema7.1}
Theorem \ref{theo1.1} and \cite{Ch4}, Theorem~2.1, show that if $(E, v)$ is a
Henselian discrete valued strictly PQL-field, then $\widehat
E(p)/\widehat E$ is a $\mathbb Z _{p}$-extension, for each $p \in
P(E)$. Therefore, one can take as $E _{\infty }$ the compositum of
all $I \in \Omega (E)$ that are inertial over $E$. Note also that if
$E$ is SQL, then Br$(E)$ is isomorphic to the direct sum $\oplus _{p
\in P(E)} \mathbb Z (p ^{\infty })$ (by (2.1) (i)), and for each $p
\in P(E)$, the set $\Lambda _{p} (E)$ can be chosen so that one may
put $C _{p} (E) = \{c _{n} (p) = \pi \colon \ n \in \mathbb N\}$,
where $\pi $ is a uniformizer of $(E, v)$. When $E$ is a local
field, and for each $(p, n) \in P(E) \times \mathbb N$, $\varphi
_{n} (p)$ is the Frobenius automorphism of the inertial extension of
$E$ in $E _{\rm sep}$ of degree $p ^{n}$, the sets $H _{p} (E)$, $p
\in P(E)$, from the proof of Theorem \ref{theo3.2}, define the Hasse (the norm
residue) symbol, in the sense of \cite{Ko2}, and give rise to the
Artin map (cf. \cite{I}, Ch. 6).
\end{rema}

\medskip
\begin{coro}
\label{coro7.2}
In the setting of (2.2) (i), let $M \in \Omega (E)$ and $F$ be a
finite extension of $E$ in $E _{\rm sep}$. Then $N(MF/F) = \{\lambda \in
F ^{\ast }\colon N _{E} ^{F} (\lambda ) \in N(M/E)\}$.
\end{coro}

\begin{proof}
It suffices to prove that $N(MF/F)$ includes the preimage of
$N(M/E)$ in $F ^{\ast }$ under $N _{E} ^{F}$. In view of Theorem
\ref{theo3.1}, Lemma \ref{lemm2.1}, \cite{Ch6}, I, Lemma~4.2 (ii), and of the
PQL-property of $F$, one may consider only the case where $M/F _{0}$
is cyclic, for $F _{0} = M \cap F$. Put $\mu _{0} = N _{F _{0}} ^{F}
(\mu )$, for each $\mu \in F ^{\ast }$. Theorem \ref{theo3.2} (iii),
\cite{Ch6}, I, Lemma~4.2 (ii), and norm transitivity in towers of
intermediate fields of $MF/E$ imply that if $N _{E} ^{F} (\mu ) \in
N(M/E)$, then $\mu _{0} \in N(M/F _{0})$. At the same time, it
follows from Proposition \ref{prop2.2} (ii), \cite{Ch6}, I, Corollary~8.5, the
RC-formula and the assumptions on $E$ that Cor$_{F/F _{0}}$ is
injective. Therefore, one deduces from the former part of (4.1)
(iii) (in the general form pointed out in Remark \ref{rema4.3}) that $\mu \in
N(MF/F)$, which proves our assertion.
\end{proof}

\medskip
Corollary \ref{coro7.2} generalizes \cite{I}, Theorem~7.6. Applying the
RC-formula, Propositions \ref{prop2.2} (ii) and \ref{prop2.3} (i), as 
well as Lemma \ref{lemm2.1}, norm and corestriction transitivity, and 
already used known relations between norms and cyclic algebras, one 
obtains by the method of proving Theorem \ref{theo3.2} the existence 
of an exact analogue to the local Hasse symbol in the following 
situation:

\medskip
\begin{coro}
\label{coro7.3}
Assume that $E$ is a nonreal field, such that every $L \in \Omega
(E)$ is strictly {\rm PQL} with $\rho _{E/L}$ surjective. Then the
maps {\rm Cor}$_{L/E}\colon \ L \in \Omega (E)$, are bijective, and
there are sets  $H(E ^{\prime }) = \{( \ , M ^{\prime }/E ^{\prime
})\colon $ $E ^{\prime \ast } \to \mathcal{G}(M ^{\prime }/E
^{\prime }) ^{{\rm Br}(E')},$ $M ^{\prime } \in \Omega (E ^{\prime
})\}$, $E ^{\prime } \in \Omega (E)$, of group homomorphisms
satisfying the following:
\par
{\rm (i)} $( \ , M ^{\prime }/E ^{\prime })$ is surjective and its
kernel equals $N(M ^{\prime }/E ^{\prime })$, for each $E ^{\prime }
\in \Omega (E)$, $M ^{\prime } \in \Omega (E ^{\prime })$;
\par
{\rm (ii)} $H _{E'}$ has the properties required by Theorem \ref{theo1.3}
(ii), for every $E ^{\prime } \in \Omega (E)$; furthermore, if $M
\in \Omega (E)$ and $K$ is an intermediate field of $M/E$, then $(\lambda
, M/K) = (N _{E} ^{K} (\lambda ), M/E)$, for any $\lambda \in K ^{\ast }$;
\par
{\rm (iii)} The sets $H(E ^{\prime }), E ^{\prime } \in \Omega (E)$,
are determined by the mappings $( \ , \Gamma /E)$, when $\Gamma $
ranges over finite extensions of $E$ in $E _{\infty }$ of primary
degrees.
\end{coro}

\medskip
Corollary \ref{coro7.3} has a partial analogue for a formally real strictly
PQL-field $E$ and the sets $\Omega ^{\prime } (E ^{\prime }) = \{M
^{\prime } \in \Omega (E ^{\prime })\colon \ 2 \dagger [M ^{\prime
}\colon E ^{\prime }]\}$, $E ^{\prime } \in \Omega ^{\prime } (E)$.
Specifically, statements (i)$\div $(iii) hold when every $E ^{\prime
} \in \Omega ^{\prime } (E)$ is $p$-quasilocal with Br$(E ^{\prime
}) _{p}$ included in the image of $\rho _{E/E'}$, for each $p \in
\mathbb P \setminus \{2\}$; also, in this case, Cor$_{E'/E}$ induces 
isomorphisms Br$(E ^{\prime }) _{p} \cong {\rm Br}(E) _{p}$, $p >
2$. This is proved in the same way as Corollary \ref{coro7.3}.  On the other
hand, it follows from the Artin-Schreier theory that if $E ^{\prime
} \in \Omega ^{\prime } (E)$ and $E ^{\prime } \neq E$, then $E
^{\prime }$ is formally real and the action of $\mathcal{G}(E
^{\prime }/E)$ induces $[E ^{\prime }\colon E]$ orderings on $E
^{\prime }$. Therefore, $E ^{\prime }$ is not $2$-quasilocal, so it
cannot admit LCFT.

\medskip
\begin{rema}
\label{rema7.4}
Under the hypotheses of Theorem \ref{theo1.2}, let $E _{\rm ab}$ be the
compositum of all $M \in \Omega (E)$, $E _{\rm ab} (p) = E _{\rm ab}
\cap E (p)$, for each $p \in P(E)$, and $N _{1} (E) = \cap _{M \in
\Omega (E)} N(M/E)$. Suppose that $d(p) \in \mathbb N$, for all $p
\in P(E)$, and Pr$_{c} (E ^{\ast })$ is the profinite completion of
$E ^{\ast }/N _{1} (E)$ (concerning its existence, see \cite{Ko1},
Sect. 1.2). Using Theorems \ref{theo1.1}, \ref{theo1.3} and Galois theory, and arguing
as in the proof of implication (iii)$\to $(i) of \cite{Ko1},
Proposition~1.14, one obtains that ${\rm Pr} _{c} (E ^{\ast })$ is
isomorphic to the topological group product $\prod _{p \in P(E)}
\mathcal{G}(E _{\rm ab} (p)/E) ^{d(p)}$, and so generalizes a part
of \cite{I}, Proposition~6.3.
\end{rema}

\medskip
Note finally that every abelian torsion group $T$ admissible by
Proposition \ref{prop2.2} (ii) is isomorphic to Br$(E(T))$, for some 
strictly PQL-field $E(T)$ satisfying the conditions of Corollary 
\ref{coro7.3} or its analogue in the formally real case. This follows 
from \cite{Ch10}, Corollary~6.6 and (7.2) (when $T$ is divisible - 
from Remark \ref{rema2.4} as well). Hence, all forms of the local 
reciprocity law and Hasse symbol admissible by Theorems 
\ref{theo1.2}, \ref{theo1.3} and Corollary \ref{coro7.3} can be 
realized.

\end{document}